\documentclass[ejs]{imsart}

\usepackage{amsmath, amsfonts, amssymb, amsthm, stmaryrd,bbm}
\usepackage{enumitem}
\usepackage[lofdepth,lotdepth]{subfig}
\usepackage{tikz}
\usepackage{pgfplots}
\usetikzlibrary{arrows}
\usepackage{algorithm2e}
\usepackage{amssymb}
\usepackage{hyperref}
\hypersetup{
    bookmarks=true,         
    colorlinks=true,       
    linkcolor=red,          
    citecolor=green,        
    filecolor=magenta,      
    urlcolor=cyan           
}

\renewcommand{\P}{\mathbb{P}}

\newlist{npf}{enumerate}{5}
\setlist[npf,1]{label={\textbf{Step \arabic*.}}, ref={Step \arabic*},leftmargin=*, itemindent=3.5em}
\newtheorem*{Hypopf}{Hypotheses in this proof}

\RequirePackage{changepage}

\newcommand{\pf}{\textsc{Proof}}

\newcommand{\onepf}[1]{\linebreak[0]\pf : #1 }

\newcommand{\sign}{\text{sign}}
\newcommand{\IF}{\mathrm{IF}}

\theoremstyle{plain}
\newtheorem{Theorem}{Theorem}
\newtheorem{Proposition}{Proposition}
\newtheorem{Corollary}{Corollary}

\newtheorem{Lemma}{Lemma}

\newtheorem{Definition}{Definition}
\newtheorem{Assumptions}{Assumption}

\renewcommand{\P}{\mathbb{P}}
\newcommand{\E}{\mathbb{E}}
\newcommand{\R}{\mathbb{R}}
\newcommand{\N}{\mathbb{N}}

\newcommand{\1}{\mathbbm{1}}

\newcommand{\argmin}{\arg\min}

\begin{document}

\begin{frontmatter}

\title{Concentration study of M-estimators using the influence function.}
\runtitle{Concentration of M-estimators.}

\author[A]{\fnms{Timothée} \snm{Mathieu}\ead[label=e1]{timothee.mathieu@inria.fr}},
\address[A]{INRIA, Scool team. Univ. Lille, CRIStAL, CNRS.\footnote{This work was done while the author was in Laboratoire Mathématiques d'Orsay for his Ph.D.} \printead{e1}}

\begin{abstract}
We present a new finite-sample analysis of M-estimators of locations in a Hilbert space using the tool of the influence function. In particular, we show that the deviations of an M-estimator can be controlled thanks to its influence function (or its score function) and then, we use concentration inequality on M-estimators to investigate the robust estimation of the mean in high dimension in a corrupted setting (adversarial corruption setting) for bounded and unbounded score functions. For a sample of size $n$ and covariance matrix $\Sigma$, we attain the minimax speed $\sqrt{Tr(\Sigma)/n}+\sqrt{\|\Sigma\|_{op}\log(1/\delta)/n}$ with probability larger than $1-\delta$ in a heavy-tailed setting. One of the major advantages of our approach compared to others recently proposed is that our estimator is tractable and fast to compute even in very high dimension with a complexity of $O(nd\log(Tr(\Sigma)))$ where $n$ is the sample size and $\Sigma$ is the covariance matrix of the inliers and in the code that we make available for this article is tested to be very fast.\end{abstract}

\begin{keyword}[class=MSC2020]
\kwd[Primary ]{62F35}
\kwd[; secondary ]{60G25}
\end{keyword}

\begin{keyword}
\kwd{Robust Statistics}
\kwd{Concentration inequalities}
\kwd{Mean estimation}
\end{keyword}

\end{frontmatter}


\section{Introduction}

Recently, there have been an increase in the quantity and dimensionality of data one has to investigate in Machine Learning tasks. Big datasets are rather difficult to handle because it is not possible anymore to check by hand some gigabyte or terabyte of data to search for possibly abnormal data, and the task is even more difficult in high dimension. An answer to this problem is robust statistics and in particular, because mean estimators are used everywhere in Machine Learning (for empirical risk minimization, cross validation, data pre-processing, feature engineering...), it has become critical to understand in depth robust mean estimators.

In this article, we try to understand one type of mean estimators: M-estimators. We show that they are nearly optimal in the heavy-tailed setting and we exhibit an algorithm that allows us to compute these estimators in a linear number of steps. Furthermore, the algorithm is very fast in practice. See Section~\ref{sec:illustration} for more illustrations of our algorithms.

The analysis of M-estimators to estimate the mean is divided in two steps: first we analyse the deviations of M-estimators around $T(P)$, which is some location parameter meant to exhibit a central tendency of the data and then we bound the bias $\|T(P)-\E[X]\|$. Let $X \sim P$ for some $P$ probability on a Hilbert space $\mathcal{H}$, let $\rho$ be an increasing function from $\R_+$ to
$\R_+$,  we are interested in estimating the location parameter $T(P)$ defined by
\begin{equation}\label{def:T_min}
T(P) \in \argmin_{\theta \in \mathcal{H}}\E\left[\rho(\|X-\theta\|)\right]=0,
\end{equation}
where $\|\cdot\|$ is a norm associated to a scalar product. Alternatively, if $\rho$ is smooth enough (which will be the case in this article), we define $T(P)$ by
\begin{equation}\label{def:T_psi}
\E\left[\frac{X-T(P)}{\|X-T(P)\|}\psi(\|X-T(P)\|)\right]=0,
\end{equation}
where $\psi=\rho'$ is called the score function. The empirical estimator obtained by plugging the empirical density $\widehat P_n$ in equation~\eqref{def:T_psi} is called M-estimator associated with $\psi$, it is denoted $T(X_1^n)$ and computed from a sample $X_1,\dots,X_n$ using the following equation:
\begin{equation}\label{def:T_hat}
\sum_{i=1}^n\frac{X_i-T(X_1^n)}{\|X_i-T(X_1^n)\|}\psi(\|X_i-T(X_1^n)\|)=0.
\end{equation}
A particular case of $T(P)$ is obtained when choosing $\psi(x)=x$ in which case $T(P)=\E[X]$ and $T(X_1^n)=\frac{1}{n}\sum_{i=1}^n X_i$, however it is well known that the empirical mean is not robust. A careful choice of the function $\psi$ yield estimators that are more robust to outliers and to heavy-tailed data (see~\cite{catoni2012}).

The subsequent problem is to see how the properties of $\psi$ impact the robustness and efficiency of $T(X_1^n)$ when estimating $T(P)$ or $\E[X]$. To study the robustness of $T(X_1^n)$ we use the influence function, a common  tool in robust statistics.
The influence function is used to quantify the robustness of an estimator, see for example \cite{hampel1974,HampelEtal86,robuststat,RONCHETTI199759} in which are derived properties such as the asymptotic variance or the breakdown point of the estimator $T(X_1^n)$ using the influence function.
The influence function is the Gâteaux derivative of $T$ evaluated in the Dirac distribution in a point $x\in \mathcal{H}$ and in the case of M-estimators in $\R^d$, from \cite[Eq~4.2.9 in Section 4.2C.]{HampelEtal86}, the influence function takes the following simple form:
\begin{equation}\label{eq:def_IF}
\IF(x,T,P)=\lim_{t \to 0}\frac{T((1-t)P+t\delta_x)-T(P)}{t}=M_{P,T}^{-1}\frac{x-T(P)}{\|x-T(P)\|}\psi(\|x-T(P)\|),
\end{equation}
where $M_{P,T}$ is a non-singular matrix whose explicit formula is not important for most of our application because of our choice of $\psi$ function (an explicit formula can however be found in the proof of Theorem~\ref{th:corruption_bias}).

The general idea is that, if the estimator is smooth enough (for example if it is Fréchet or Hadamard differentiable, see~\cite{fernholtz}), then one can write the following expansion
\begin{equation}\label{eq:taylor_informal}
T(P)=T(Q)+\int_{\mathcal{H}} \IF(x,T,Q)d(P-Q)(x)+R(P,Q),
\end{equation}
where the remainder term $R(P,Q)$ is controlled.
For example, if we apply equation~\eqref{eq:taylor_informal} to $Q=\widehat P_n$ the empirical distribution, the influence function provides a first order approximation for the difference between the estimator $T(\widehat P_n)=T(X_1^n)$ and its limit $T(P)$.
This technique of approximating the estimator by its influence function is also linked to the Bahadur decomposition, see~\cite{bahadur1966} and~\cite{he1996} for applications to M-estimators.
The influence function of $M$-estimators is usually chosen bounded in robust statistics, in particular from \cite{hampel1971,robuststat} we have that if $\psi$ is bounded, then the influence function is bounded and $T$ is qualitatively robust (i.e. the estimator $T(X_1^n)$ is equi-continuous, c.f. \cite{robuststat}) and have asymptotic breakdown point $1/2$. On the other hand if $\psi$ is unbounded, then $T(X_1^n)$ is not qualitatively robust, the influence function is not bounded and the asymptotic breakdown point is zero. From Hampel's Theorem~\cite[Theorem 2.21]{robuststat} we also have that $\psi$ is bounded  if and only if $T$ is a continuous functional with respect to the Levy metric. More generally, the influence function has been used in a lot of works on asymptotic robustness, see~\cite{HampelEtal86,robuststat} or~\cite{hampel1974,RONCHETTI199759}.

The influence function has also been used recently in Machine Learning literature in order to have a model selection tool specialized in robustness, see for example~\cite{debruyne2008model}, \cite{koh2017understanding} and the closely related tool of leave one out error~\cite{looarticle}.
The field of Robustness in Machine Learning has been very active in the last few years, in particular after several works by Olivier Catoni and co-authors  in~\cite{catoni2012,catoni2017dimensionfree}, the goal is to prove non-asymptotic deviation bounds when the data are more heavy-tailed than what is usually considered in classical Machine Learning. This line of thought has been continued in a number of articles, in particular \cite{subgaussian} introduced some general concept of sub-Gaussian estimators that have been then used successfully in other applications, see~\cite{chinot2018statistical, 1911.05911,catoni2012, DBLP:journals/focm/LugosiM19,zhou2018,1910.07485}. See also some comprehensive lecture notes on the subject in \cite{1908.10761}.

It is interesting to note that contrary to works from classical robust theory from the 70's, the influence functions of the M-estimators used by Catoni are not necessarily bounded. In this article, we initiate the analysis of the effect of unbounded influence function on the robustness of M-estimators, Huber~\cite{huber1964} told us that the influence function must be bounded while Catoni uses unbounded influence function and he still shows robust properties for this type of estimator, the difference is in their vision of what is a robust estimator.

There will be three parts in our analysis of this problem, first we extend Catoni's non-asymptotic analysis of M-estimators as we analyze M-estimators with more general influence functions and in a multivariate setting using the properties of the influence function, making a link between the deviations of the influence function and the deviations of $T(X_1^n)$. Second we apply our theory to three specific M-estimators for which we show tight upper bound on the rate of convergence to the mean. Finally, we investigate an algorithm to compute $T(X_1^n)$ and we show that this algorithm converges in a reasonable number of steps.

More precisely, In Section~\ref{sec:tail}, we show that concentration inequalities for M-estimators derive from concentration inequalities on the influence function by showing roughly that
\begin{equation}\label{eq:main_result_sec_3}
\|T(X_1^n)-T(P)\|  \simeq  \left\|\frac{1}{n}\sum_{i=1}^n \frac{X_i-T(P)}{\|X_i-T(P)\|} \psi(\|X_i-T(P)\|)\right\| .
\end{equation}
From equation~\eqref{eq:def_IF}, the right hand side of equation \eqref{eq:main_result_sec_3} can be interpreted as the deviation of the influence function.
The right hand side of equation~\eqref{eq:main_result_sec_3} can be controlled under classical assumptions. For example in $\R$, if $\psi$ is bounded by $\beta>0$ (Huber estimator), we can use Hoeffding or Bernstein inequality to get a control on $\|T(X_1^n)-T(P)\|$. Using Hoeffding inequality, we obtain a concentration rate similar to the rate of the empirical mean on Gaussian data.
$$\P\left(\left|\frac{1}{n}\sum_{i=1}^n \psi(X_i-T(P))\right| \ge \frac{\beta}{\sqrt{n}}\lambda\right)\le e^{-2\lambda^2}.$$
Remark that using Bernstein inequality, we don't need to have $\psi$ bounded, hence this gives us a mean to show sub-Gaussian rates for M-estimators with unbounded $\psi$ function.


%

In Section~\ref{sec:bias}, we provide bounds on the bias $\|T(P)-\E[X]\|$ and on the variance terms in the concentration inequality from Section~\ref{sec:tail}. Bounding the bias has often been a problem in robust statistics, if the distribution is skewed and the bias is not controlled we can only say that we estimate a quantity meant to quantify a central tendency of $P$ but we do not estimate $\E[X]$ directly. However in statistical learning for example, estimating the mean is not just an arbitrary choice and we don't want to estimate a central tendency of the dataset, we want to estimate its mean. In this article, we give explicit bounds on the bias and we use those bounds (in Section~\ref{sec:examples_tail}) to give concentration results on $T(X_1^n)$ around $\E[X]$ in the context of heavy-tailed and adversarially corrupted datasets even beyond the $L^2$ case. Indeed, we show that the $L^2$ assumption is needed only to be able to handle the bias $\|T(P)-\E[X]\|$ but on the other hand if the distribution is symmetric, $T(P)=\E[X]$ and we can obtain sharp deviation bounds even in the case of $L^1$ distributions. Bounds on the bias are already present for the specific case Huber estimator in regression in \cite{sun2020adaptive}. Our bound on the bias has the same flavor as \cite[Proposition 1]{sun2020adaptive} but applied to multivariate location parameter for more general M-estimators. Similarly to the discussion in \cite{sun2020adaptive}, our analysis leads us to the study of a bias-variance trade-off for M-estimator depending on the value of $\beta$.

In this context, in Section~\ref{sec:examples_tail}, we show that $T(X_1^n)$ is suitable to estimate the mean in high dimension in a heavy-tailed and corrupted setting (even though our estimators are not minimax in corrupted setting). In the literature, there are estimators that have strong theoretical guarantees but that are intractable, for example one can see estimators based on the aggregation of one-dimensional estimators (same idea as projection pursuits), see \cite[Theorem 44]{1908.10761} and reference therein, see also~\cite{lugosi2019sub} and there has also been estimators based on depth, for example Tukey's median~\cite{chen2018robust}. On the other hand, there are tractable algorithms but with non minimax optimal rates of convergence for example the coordinate-wise median or the geometrical median~\cite{MR3378468,catoni2017dimensionfree}, our work belong to this type of methods, our estimator is easily computable and even though the obtained error bounds are much better than for the coordinate-wise median, at least in corrupted setting it is not minimax. 
Recently there have been several propositions of algorithms whose goal was to be at the same time tractable and minimax, see~\cite{diakonikolas2020outlier,zhu2020robust,hopkins2020robust,depersin2019robust,hopkins2020mean,cherapanamjeri2019fast} however these algorithms are often hard to implement and in practice the complexity makes them intractable for high-dimensional problems.

In a corrupted setting where inliers have a $q^{th}$ finite moment, for $q>2$, we control the deviations of M-estimators. Let us suppose that the data are an adversarial $\varepsilon_n$ corrupted dataset (see Assumption~\ref{ass:adv}). Then, under some assumptions on $\psi$,  for all $0<t\lesssim n \frac{\sqrt{\|\Sigma\|_{op}}}{ \E\left[\| X- \E[X]\|^q\right]^{1/q}}$, with probability larger than $1-8\exp(-t)$,
\begin{align*}
\left\|\E[X]-T(X_1^n)\right\|  \lesssim &   \frac{\sqrt{\mathrm{Tr}(\Sigma)}+\sqrt{\|\Sigma\|_{op}t}}{\sqrt{n}} \bigvee \E\left[\| X- \E[X]\|^q\right]^{1/q} \varepsilon_n^{1-1/q}.
\end{align*}
See Proposition~\ref{prop:huber_corrupted_q} below for the formal and more precise statement in the case of Huber's estimator, and see Section~\ref{sec:examples_tail} for other examples in particular for M-estimators with unbounded score function. In the heavy-tailed setting where $\varepsilon_n=0$, the bound is almost minimax optimal, the difference is that $t $ is allowed to be up to $n \frac{\sqrt{\|\Sigma\|_{op}}}{ \E\left[\| X- \E[X]\|^q\right]^{1/q}}$ instead of the usual $t \lesssim n$. This estimator constitutes one of the few tractable and efficient estimators of the multivariate mean in heavy-tailed setting.

The error due to corruption is of order $O(\E\left[\| X- \E[X]\|^q\right]^{1/q} \varepsilon_n^{1-1/q})$, the term $\varepsilon_n^{1-1/q}$ is optimal (see \cite[Lemma 5.4]{minsker2018uniform}) but the factor $\E\left[\| X- \E[X]\|^q\right]^{1/q}$ is sub-optimal. In Section~\ref{sec:lower_bound} we show that the dependency in the dimension due to the $\E\left[\| X- \E[X]\|^q\right]^{1/q}$ is unavoidable (to understand the magnitude of this term, one may think about the case where all the marginal have same law, then we can show that $\E\left[\| X- \E[X]\|^q\right]^{1/q}\lesssim \sqrt{d}$). We can't avoid this dimension dependence for our estimator contrary to~\cite{depersin2019robust,lai2016agnostic} in which the authors show that the optimal error due to corruption is $O(\varepsilon\sqrt{\|\Sigma\|_{op}})$ in the case of Gaussian distributions. Hence our estimator is not optimal with respect to corruption. 

Finally, in Section~\ref{sec:algo}, we exhibit an algorithm to compute $T(X_1^n)$ and we show that this algorithm converges in a finite number of steps resulting in a complexity of order $O(nd\log(\mathrm{Tr}(\Sigma)))$ where $n$ is the sample size and $\Sigma$ is the covariance matrix of the inliers (our analysis is only valid for $d \ll e^n$). This algorithm is practically efficient, it is illustrated in Section~\ref{sec:illustration} and we invite the interested reader to check out the github repository~\url{https://github.com/TimotheeMathieu/RobustMeanEstimator} for the python code. 
\section{Setting and Notations}\label{sec:setting_if}
\subsection{Setting}\label{sec:examples}
In all the article, we consider $X_1,\dots, X_n$ in a Hilbert space $\mathcal{H}$ that have been corrupted by an adversary using the following process.

\begin{Assumptions}\label{ass:adv}
There exist $X_1',\dots,X_n' \in \mathcal{H}$ i.i.d. following a law $P$ that have been modified by an ``adversary" to obtain $X_1,\dots,X_n$. The adversary can modify at most $|\mathcal{O}|$ points such that there exist $\mathcal{I},\mathcal{O}$ partition of $\{1,\dots, n\}$, with $\mathcal{I}\cup \mathcal{O} = \{1,\dots, n\}$ and $|\mathcal{O}|<|\mathcal{I}|$.
\end{Assumptions}

The points $(X_i)_{i \in \mathcal{O}}$ are arbitrary and are called outliers, the points $(X_i)_{i \in \mathcal{I}}$ are called inliers. Remark that the statistician do not know the sets $\mathcal{I}, \mathcal{O}$ and that although we have $X_i = X'_i$ for any $i \in \mathcal{I}$, the sample $(X_i)_{i \in \mathcal{I}}$ is in general not i.i.d. because the adversary can choose which data are corrupted using knowledge on the inliers. For instance it is possible that the adversary decided to corrupt the $|\mathcal{O}|$ points from $X_1', \dots,X_n'$ that were the closest to the theoretical mean in which case the $(X_i)_{i \in \mathcal{I}}$ are not independent.

We consider the functional $T$ defined by
\begin{equation}\label{eq:def_T_th}
  \E\left[\frac{X-T(P)}{\|X-T(P)\|}\psi(\|X-T(P)\|)\right]=0,
\end{equation}
for some $\psi: \R_+ \to \R_+$, the existence and unicity of $T(P)$ is discussed in Lemma~\ref{lem:unicity}. We are interested in the behavior of the associated M-estimator $T(X_1^n)$ defined by
\begin{equation}\label{eq:def_T_emp}
\sum_{i=1}^n \frac{X_i-T(X_1^n)}{\|X_i-T(X_1^n)\|}\psi(\|X_i-T(X_1^n)\|)=0,
\end{equation}
where $\psi$ is a function that satisfies the following properties.
\begin{Assumptions}\label{ass:prop_psi}
The function $\psi:\R_+ \to \R_+$ verifies the following properties:
\begin{enumerate}[label=(\roman{*}), ref=2-(\roman{*}), leftmargin=*]
\item $\psi$ is continuous and differentiable almost everywhere \label{ass:1}
\item $\psi(0)=0$ \label{ass:2}
\item $\psi$ is concave \label{ass:concave}
\item There exist $\beta,\gamma>0$ such that \label{ass:3}
$$
\forall x\ge 0,\qquad 1 \ge \psi'(x)\ge \gamma\1\{x\le \beta \},
$$
where $\1$ is the indicator function.
\end{enumerate}
\end{Assumptions}

When we want to emphasize the dependency in $\beta$, we will use the notation $\psi_1$ to be such that for all $x\ge 0$, $\psi(x)=\beta\psi_1(x/\beta)$.

Because $\psi$ is concave, non decreasing and not identically zero, there are always a couple of positive constants $\beta,\gamma$ such that Assumptions~\ref{ass:prop_psi} holds. For our results to hold we will ask that $\beta$ and $\gamma$ are not too small.
A first result that can be derived from Assumptions~\ref{ass:prop_psi} and  some additional assumptions is that our problem is well defined. This comes from the fact that the problem is a convex problem.

Let
$$Z_{P,\psi}:\theta \mapsto \E_{X\sim P}\left[\frac{X-\theta}{\|X-\theta\|} \psi\left(\left\|X-\theta\right\|\right) \right],$$
we have the following lemma whose proof is in Section~\ref{sec:proof_convexity}.
\begin{Lemma}\label{lem:convexity}
Let $\psi$ satisfy Assumptions~\ref{ass:prop_psi}, let $u \in S$ and $\theta \in \mathcal{H}$,
\begin{align*}
u^TJac(Z_{P,\psi})(\theta)u\le-\E\left[\psi'\left(\left\|X-\theta\right\|\right) \right]
\end{align*}
where $Jac$ denotes the Jacobi operator.
\end{Lemma}
Using the previous lemma, we can prove using additional hypotheses that the problem is well defined.

\begin{Lemma}\label{lem:unicity}
Let $\psi$ satisfy Assumptions~\ref{ass:prop_psi}, define $\rho :x \mapsto \int_{0}^x \psi(t)dt$ and let $X$ satisfy $\E[\rho(\|X-\E[X]\|)]<\rho(\beta)$, then $T(P)$ defined by equation~\eqref{eq:def_T_th} exists and is unique.
\end{Lemma}
This lemma is proven in Section~\ref{sec:proof_unicity}. Remark that the condition $\E[\rho(\|X-\E[X])\|)]<\rho(\beta)$ is implied by $\mathrm{Tr}(\Sigma)=\E[\|X-\E[X]\|^2]< 2\rho(\beta)$ because $\rho(x)\le x^2/2$ however an hypothesis on $\mathrm{Tr}(\Sigma)$ is a lot stronger because it supposes a finite second moment.
In the whole article, we will suppose that $T(P)$ is unique, we do not necessarily suppose that the assumptions of Lemma~\ref{lem:unicity} hold because they are not minimal assumptions for unicity and existence of $T(P)$. However, for simplicity, we will suppose that the following condition is verified.

\begin{Assumptions}\label{ass:unicity}
$X$ satisfies $\E[\rho(\|X-\E[X]\|)]< \min(\rho(\beta/3), \psi(\beta/2)^2/2)$.
\end{Assumptions}
Assumption~\ref{ass:unicity} means that $\beta$ should be large enough to encompass most of the weight of $X$ around the expectation. When $\psi$ is bounded, this is only a first moment assumption. This assumption is due to limitations in our methods and we believe that this is not necessary. One can show that $\rho(x)\ge \psi(x)^2/2$ which may be used to simplify Assumption~\ref{ass:unicity}.

Assumptions~\ref{ass:prop_psi} and Assumptions~\ref{ass:unicity} will be supposed true. The behavior of $\psi$ at $0$ allows us to control the deviations of the estimator using the influence function, see Section~\ref{sec:tail} and it is also important to control the bias of the resulting estimator, see Section~\ref{sec:bias}.
On the other hand, the growth rate of $\psi$ at $+\infty$ is central to derive concentration bounds of $T(X_1^n)$, as will become clear all along Section~\ref{sec:tail} and~\ref{sec:examples_tail}.
Assumptions~\ref{ass:prop_psi} do not always apply to M-estimators, for example the sample median is not an estimator derived from a function $\psi$ satisfying these assumptions.
On the other hand, we provide three examples of score functions satisfying Assumptions~\ref{ass:prop_psi}, with three different growth rates when $x$ goes to infinity.

\begin{description}

\item[Huber's estimator.] Let $\beta>0$. For all $x\ge 0$, let
\begin{equation}\label{eq:huber_score_1}
\psi_H(x)=
x \,\1\{x\le \beta\}+\beta \, \1\{x>\beta\}.
\end{equation}
In dimension $1$, the M-estimator constructed from this score function is called Huber's estimator~\cite{huber1964}.

\item[Catoni's estimator.] Let $\beta>0$. For all $x\ge 0$, let
\begin{equation}\label{eq:catoni_score}
\psi_C(x)=\beta\log\left(1+\frac{x}{\beta}+\frac{1}{2}\left(\frac{x}{\beta}\right)^2\right).
\end{equation}
The associated M-estimator is one of the estimators considered by Catoni in~\cite{catoni2012}.
We call the resulting M-estimator Catoni's estimator.

\item[Polynomial estimator.] Let $p\in \N^*$, $\beta>0$. For all $x\ge 0$, let
\begin{equation}\label{eq:poly_score}
\psi_{P}(x)=\frac{x}{1+\left(\frac{x}{\beta}\right)^{1-1/p}}.
\end{equation}
We call Polynomial estimator the M-estimator obtained using this score function.
\end{description}

\begin{figure}[h]
\begin{center}
\includegraphics[scale=0.66]{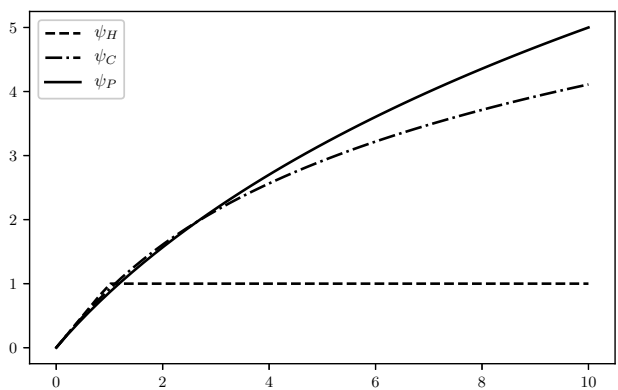}
\caption{Plot of $\psi_H$ and $\psi_C$ for $\beta=1$. $\psi_P$ is plotted for $\beta=10$ and $p=5$.}
\end{center}
\end{figure}
The following result shows that the score functions from the previous three examples satisfy Assumptions~\ref{ass:prop_psi}.
\begin{Lemma}\label{lem:ass_exemples}
For all $x\ge 0$, we have
$\psi'_H(x)= \1\{x\le \beta\}, $
$\psi'_C(x)\ge \frac{4}{5}\1\{x\le \beta\}, $
$\psi_P'(x)\ge \frac{1}{4}\1\{x\le \beta\}. $
\end{Lemma}
The proof of Lemma~\ref{lem:ass_exemples} is postponed to Section~\ref{sec:proof_ass_exemples}.

\subsection{Notations}\label{sec:notations}
Let $\mathcal{P}$ denote the set of probability distributions on $\mathcal{H}$, $S=\{x\in \mathcal{H}:\|x\|=1\}$ where $\|\cdot\|$ is a norm associated with a scalar product $\langle \cdot,\cdot\rangle$ on $\mathcal{H}$. For any $\psi:\R_+\to\R_+$, let $\mathcal{P}_\psi=\{ P\in \mathcal{P}: \E_P[\psi(\|X\|)]<\infty\}$. We denote $a\lesssim b$ if there exists a numerical constant $C>0$ such that $a \le C b$.

Let $T_H$, $T_C$ and $T_P$  denote the functionals associated to the score functions $\psi_P$,$\psi_H$ and $\psi_C$ respectively, using Equation~\eqref{eq:def_T_th}. Define the following variance terms
$$V_\psi=\E[\psi(\|X-T_H(P)\|)^2], \quad \text{and}$$
$$v_\psi=\left\|\E\left[\frac{(X-T_H(P))(X-T_H(P))^T}{\|X-T_H(P)\|^2}\psi(\|X-T_H(P)\|)^2\right]\right\|_{op}.$$
 These variance terms are to be compared with $Tr(\Sigma)$ and $\|\Sigma\|_{op}$ in the multivariate Gaussian setting, Hanson-Wright inequality (see Equation~\ref{eq:hanson_wright}) tells us that $Tr(\Sigma)$ and $\|\Sigma\|_{op}$ describe the spread of the empirical mean in high dimension. Here we are not in a Gaussian setting and for example in the case of Huber's estimator, $V_{\psi_H}$ and $v_{\psi_H}$ will describe the spread of the influence function of Huber's estimator. See Section~\ref{sec:bound_variance} for a more formal study of the link between $V_{\psi_H},v_{\psi_H}$ and $Tr(\Sigma),\|\Sigma\|_{op}$.

\section{Tail probabilities of M-estimator and Influence function}\label{sec:tail}
\label{sec:tail_general}

The main result of the paper compares the tail probabilities of $\|T(X_1^n)-T(P)\|$ and those of its influence function.

\begin{Definition}\label{def:tail_fun}
We call $t_{T}$ and $t_{\IF}$ the tail probability functions defined for all $\lambda>0$ and $\theta \in \mathcal{H}$ by
\begin{gather}
\label{eq:t_T} t_T(\lambda; \theta, X_1^n):=\P\left(\|T(X_1^n)-\theta\|\ge\lambda \right)\\
t_{IF}(\lambda;\theta, X_1^n):=\P\left(\left\|\frac{1}{n}\sum_{i=1}^n \frac{X_i-\theta}{\|X_i-\theta\|} \psi(\|X_i-\theta\|)\right\| \ge\lambda\right).
\end{gather}
\end{Definition}
The main theorem of Section~\ref{sec:tail} is the following.
\begin{Theorem}\label{th:tIF_tT_1}
Suppose that Assumptions\ref{ass:adv}, \ref{ass:prop_psi} and \ref{ass:unicity} hold. If moreover $V_\theta=\E_P[\psi(\|X'-\theta\|)^2]\le \psi(\beta/2)^2/2<\infty$ and $|\mathcal{O}|\le n\gamma/8$, then for all $\lambda \in (0,\beta/2)$ and for all $\theta \in \mathcal{H}$,
\begin{equation}\label{eq:main_result_concentration1}
t_T(\lambda; \theta, X_1^n)\le t_{\IF}\left(\lambda\gamma/4; \theta, X_1^n\right)+e^{-n\gamma^2/32}.
\end{equation}
\end{Theorem}
The proof of this result is given in Section~\ref{sec:proof_tIF_tT}. In the heavy-tailed setting, we will use Theorem~\ref{th:tIF_tT_1} with $\theta=T(P)$ in order to have a small value of $t_{\IF}(\lambda;\theta,X_1^n)$. On the other hand, in a corrupted setting, $\theta$ will be set to $T(P)$ where $P$ is the law of inliers. For now, there is no hypothesis on the outliers in $\mathcal{O}$, in what follows we will see that if $\psi$ is unbounded, we will need some hypothesis on $X_i, i\in \mathcal{O}$ in order to have a control of the value of $t_{\IF}$.

Remark that although M-estimators with bounded $\psi$ are proven to have a breakdown point of $1/2$ (see \cite{robuststat}), our result is only valid for a proportion of outliers $|\mathcal{O}_n|/n\le \gamma/8$. This is an artifact of the proof and with more stringent condition on $\beta$, we could allow for a higher breakdown point at the cost of a more complicated analysis. Remark however that when the corruption is close to $1/2$, the error due to corruption is the principal source of error and a looser deviation analysis could then be allowed to reach a higher breakdown point.

\textbf{Example:} Huber's estimation in dimension $d=1$. In the case of Huber's estimator, $\psi$ is bounded and from Lemma~\ref{lem:ass_exemples}, $\gamma=1$. Because $\psi_H$ is bounded by $\beta$, we have directly from Bernstein inequality for all $t>0$,
$$\P\left(\left| \frac{1}{n}\sum_{i=1}^n \sign(X-T_H(P))\psi_H(|X-T_H(P)|)\right|> \sqrt{\frac{2V_{\psi_H} t}{n}}+\frac{\beta t}{n}\right) \le 2e^{-t}.$$
Then, by Theorem~\ref{th:tIF_tT_1}, if $V_{\psi_H}\le \beta^2/8$, for all $t>0$ such that $4\sqrt{2V_{\psi_H} t/n}+4\beta t/n \le \beta/2$,
\begin{equation}\label{eq:huber_1D_1}
\P\left(\left| T_H(X_1^n)-T_H(P)\right|> 4\sqrt{\frac{2V_{\psi_H} t}{n}}+4\frac{\beta t}{n}\right) \le 2e^{-t}+e^{-n/32}.
\end{equation}
Remark that choosing $\beta=\sqrt{V_{\psi_H}n}$ gives us a sub-Gaussian concentration around $T_H(P)$, this is similar to the concentrations inequalities introduced in~\cite{subgaussian} except that we concentrate around $T_H(P)$ instead of $\E[X]$. Remark also that the condition $V_{\psi_H}\le \beta^2/8$ is rather weak because we already have $V_{\psi_H}\le \beta^2$, the condition asks that there is enough weight in the interval $[-\beta,\beta]$.
%
%

\section{Bias and variance of M-estimators when considered as estimators of the mean}\label{sec:bias}
If $P$ is symmetric then we can avoid the problem of the bias and $T(P)=\E[X]$, but unfortunately in the case of skewed distribution the bias $\|T(P)-\E[X]\|$ can be very large and the choice of $\beta$ will determine how large the bias is. In this section, we show how the bias behaves as $\beta$ grows, we also provide bounds on the variance terms defined in Section~\ref{sec:notations}, these bounds will be useful to derive concentration inequalities on $T(X_1^n)$. We will use the notation $\psi(x)=\beta\psi_1(x/\beta)$ when we want to emphasize the dependency on $\beta$.

We introduce the following function
$$Z_\beta:\theta \mapsto \E\left[\beta\frac{(X-\theta)}{\|X-\theta\|} \psi_1\left(\frac{\left\|X-\theta\right\|}{\beta}\right) \right] .$$
The following theorem links $Z_\beta$ with the distance between $T(P)$ and $\E[X]$.
\begin{Theorem}\label{th:cornerstone} 
Let $X$ be a random vector in $\mathcal{H}$, $X\sim P$ with finite expectation and suppose Assumptions~\ref{ass:prop_psi} and \ref{ass:unicity}. Then we have $ \|\E[X]-T(P)\|\le \frac{2}{\gamma}\left\|Z_\beta(\E[X])\right\|.$
\end{Theorem}
We postpone the proof to Section~\ref{sec:proof_cornerstone}. From Theorem~\ref{th:cornerstone}, it is sufficient to upper bound $\|Z_\beta(\E[X])\|$ to get a bound on the bias.



\subsection{Bias of M-estimators}
We begin with the bias of Huber's estimator obtained from equation~\eqref{def:T_hat} with $\psi=\psi_H$. The following lemma gives a bound on the bias of Huber's estimator for a distribution with a finite number of finite moments.

\begin{Lemma}\label{lem:bias_huber_using_concentration}
Let $X$ be a random variable with $\E[\|X\|^q]<\infty$ for $q\in \N^*$ and suppose that $\rho_1(1/3)\ge \E[\rho_1(\|X-\E[X]\|/\beta)]$. Then
$$\|\E[X]-T_H(P)\|\le  \frac{2\E[\|X-\E[X]\|^q]}{(q-1)\beta^{q-1}}.$$
\end{Lemma}
The proof is in Section~\ref{sec:proof_bias_huber}. Lemma~\ref{lem:bias_huber_using_concentration} is not exactly tight as can be seen for instance if we do the computation with the Pareto distribution for $d=1$ for which if the shape parameter is $\alpha=2$, we have only one finite moment but can show that in this parametric case, Huber's estimator achieves rates $1/\beta$.

The choice of $\beta$ is a very important problem when estimating $\E[X]$ using $T(X_1^n)$ and in particular we will need to choose $\beta$ carefully as a function of $n$ in order to have $T(X_1^n)$ that converges to $\E[X]$. The choice of $\beta$ will entail a sort of bias-variance tradeoff. Remark that we do not need a finite second moment for our analysis to work, we only need $\E[\rho(\|X-T(P)\|/\beta)]<\infty$ which translates in a finite first moment in the case of $\psi=\psi_H$.

In addition to Lemma~\ref{lem:bias_huber_using_concentration} we can also show an exponential bound on the bias when the random variable $X$ is sub-exponential however because the primary use of Huber's estimator is with robust statistics, we only state the result for a finite number of finite moments as it is what will interest us. An interested reader can adapt the proof to lighter-tailed distributions.

For a $\psi$ function that is not Huber's score function, the bias also depends strongly on the behavior of $\psi$ near $0$.
\begin{Lemma}\label{lem:bias_general}
Suppose that $\psi$ is $k$ times differentiable with bounded $k^{th}$ derivative and that Assumptions~\ref{ass:prop_psi} and  \ref{ass:unicity} hold, $\psi'(0)=1$ and for $2\le j\le k-1$, $\psi^{(j)}(0)=0$. Let $X$ be a random variable such that $\E[\|X\|^{k}]<\infty$, then,
\begin{equation}
\|T(P)-\E[X]\|\le \frac{\|\psi_1^{(k)}\|_\infty}{\gamma k!\beta^{k-1}}\E\left[\|X-\E[X]\|^{k} \right].
\end{equation}
Moreover, if $X$ follows a Bernoulli distribution of parameter $p$, this bound is tight in its dependency in $\beta$. When $\beta \to \infty$, we have
$$Z_\beta(\E[X])=\psi_1^{(k)}(0)\frac{p(1-p)^k-(1-p)p^k}{k!\beta^{k-1}}+o\left(\frac{1}{\beta^{k-1} }\right) $$
\end{Lemma}

This Lemma is proven in Section~\ref{sec:proof_bias_general} . For example, we can show that for Catoni's score function $\psi(x)=\log(1+x+x^2/2)$ whose second derivative is $\psi''(x)=-(x+x^2/2)/(1+x+x^2/2)^2$, we have that $\psi(x)=x-x^3/6+o(x^3)$ and then the bias of Catoni's estimator is in general of order $1/\beta^2$. Lemma~\ref{lem:bias_general} shows that the bias depends on the smoothness of the function near $0$ and also the number of finite moments.

\subsection{Bound on the variance of M-estimators}\label{sec:bound_variance}
First, we have to control the variability of $T(X_1^n)$ in order to control its deviations. The following lemma gives an upper bound on both $V_\psi$ and $v_\psi$ defined in Section~\ref{sec:notations}.
\begin{Lemma}\label{lem:bound_variance_variance}
Suppose that Assumptions~\ref{ass:prop_psi} and~\ref{ass:unicity} are satisfied, suppose that $X$ has a finite second moment with covariance operator $\Sigma$ we have that
$V_\psi\le \E[\|X-\E[X]\|^2]=\mathrm{Tr}(\Sigma),$ and $v_\psi\le \|\Sigma\|_{op}+\|\E[X]-T(P)\|^2$.
\end{Lemma}
Lemma~\ref{lem:bound_variance_variance} (proven in Section~\ref{sec:proof_bound_variance}) gives a control on $V_\psi$ and $v_\psi$ using the properties of $X$. Next we show that Lemma~\ref{lem:bound_variance_variance} is tight in the case of Huber's estimator as long as $X$ is sufficiently concentrated using the following lemma whose proof is provided in Section~\ref{sec:proof_lower_variance}.

\begin{Lemma}\label{lem:lower_bound_variance_huber}
Suppose that Assumptions~\ref{ass:prop_psi} and \ref{ass:unicity} are satisfied and that $X$ is such that $\E[\|X\|^{2q}]<\infty$ for some $q>1$, then
\begin{align*}
V_{\psi_H}\ge \E[\|X-\E[X]\|^2]-4^{q}\frac{\E\left[\|X-T_H(P)\|^{2q}\right]^{1-1/q}}{\left(\E\left[\|X-T_H(P)\|^{2q}\right]+\beta^{2q}\right)^{1-2/q}}.
\end{align*}
and
$$V_{\psi_H}\ge \|\Sigma\|_{op} -4^{q}\frac{\E\left[\|X-T_H(P)\|^{2q}\right]^{1-1/q}}{\left(\E\left[\|X-T_H(P)\|^{2q}\right]+\beta^{2q}\right)^{1-2/q}} $$
\end{Lemma}
Lemmas~\ref{lem:bound_variance_variance} and~\ref{lem:lower_bound_variance_huber} imply that if $X$ has enough moments, say with $4$ finite moments, and if $\beta$ is sufficiently large, then the behavior of the variance term is the same as the variance term for the empirical mean. On the other hand, if $X$ is not very concentrated, Lemma~\ref{lem:bound_variance_variance} can be a very rough bound and in the case of Huber estimator if $X$ has only a finite first moment but no finite variance, then $V_{\psi_H}$ and $V_{\psi_H}$ are finite even though $Tr(\Sigma)=\|\Sigma\|_{op}=\infty$.

\section{Application to the concentration of M-estimators around the mean in Corrupted Datasets}\label{sec:examples_tail}
In this section, we investigate the concentration of the three M-estimators taken as example in this article in a corrupted, heavy-tailed setting. The goal will be to recover deviations similar to the one we would have in a Gaussian setting, but when the data are not Gaussian. The gold standard in this context is the deviation of the empirical mean in a Gaussian setting (see~\cite{concentration}). If $X_1,\dots,X_n$ are i.i.d from $\mathcal{N}(\mu,\sigma^2)$ for some $\mu \in \R$ and $\sigma>0$, then for all $t>0$,
\begin{equation}\label{eq:gaussian_dev}
\P\left(\left|\frac{1}{n}\sum_{i=1}^n X_i-\mu \right|>\sigma\sqrt{\frac{t}{2n}} \right)\le e^{-t}.
\end{equation}
An equivalent of this in the multi-dimensional setting is Hanson-Wright inequality~\cite{hanson1971}: let $X \sim \mathcal{N}(\mu,\Sigma)$ for $\Sigma$ a positive definite matrix, $\mu\in \R^d$. Then, for any $t>0$,
\begin{equation}\label{eq:hanson_wright}
\P\left(\left\|\frac{1}{n} \sum_{i=1}^n X_i-\mu\right\|^2 > \frac{2Tr(\Sigma)}{n}+\frac{9t\|\Sigma\|_{op}}{n} \right)\le e^{-t}.
\end{equation}
This form of Hanson-Wright inequality can be found for example in~\cite{1908.10761}. Our aim is to obtain deviations similar to the ones in equations~\eqref{eq:gaussian_dev} and \eqref{eq:hanson_wright} but in a non-Gaussian setting.

The results we show in this section are not optimal, they are nearly optimal in the heavy-tailed setting but the effect due to corruption is sub-optimal compared to~\cite{depersin2019robust}. One of our goals is to illustrate the use of the influence function and particularly Theorem~\ref{th:tIF_tT_1} for an easy derivation of concentration inequalities for M-estimators.
We also illustrate an interesting phenomenon derived from Theorem~\ref{th:tIF_tT_1} by showing that the concentration of $T(X_1^n)$ around $T(P)$ can be much faster than the concentration of $T(X_1^n)$ around $\E[X]$, because the variance term is not $Tr(\Sigma)$ but $V$ which can be a lot smaller than $Tr(\Sigma)$ (for instance if $Tr(\Sigma)$ is not finite). It also shows that $T(X_1^n)$ is never arbitrarily bad as long as $\E[\psi(\|X\|)]<\infty$ contrary to the empirical mean.
We use the following corollary of \cite[ Theorem 4]{adamczak2008} recalled in Section~\ref{sec:reminder}, this lemma is proven in Section~\ref{sec:proof_subgaussian}.

\begin{Lemma}\label{lem:subgaussian}
Let $Y_1,\dots,Y_n$ be i.i.d random variables taking values in $\mathcal{H}$, centered and with covariance operator $\Sigma$, and such that the Orlicz norm of $Y$ is finite:
$$\|Y\|_{\psi_1}=\inf\left\{\lambda >0:\, \E\left[\exp\left(\|Y\|/\lambda \right)-1\right]\le 1 \right\} <\infty.$$
There exists an universal constant $C>0$ such that, for all $t\ge 0$,
\begin{equation}\label{eq:concentration}
\P\left(\left\|\sum_{i=1}^n Y_i\right\|\ge \frac{3}{2}\sqrt{\E\left[\left\|Y\right\|^2\right]n}+2\sqrt{nt\|\Sigma\|_{op}}+Ct\| \max_{1\le i\le n} \|Y_i\| \|_{\psi_1} \right)\le 4\exp\left(-t \right).
\end{equation}
\end{Lemma}

The last term in equation~\eqref{eq:concentration} can be handled using \cite[Lemma 2.2.2]{MR1385671} from which we get that there exists an absolute constant $K>0$ such that
\begin{equation}\label{eq:psi1_max}
\left\| \max_{1\le i\le n} \|Y_i\| \right\|_{\psi_1}\le K\log(n) \left\|Y_i \right\|_{\psi_1}.
\end{equation}
However, note that Hanson-Wright's inequality for Gaussian random variables shows that this logarithm factor is not optimal. This extra logarithm factor can be removed if $Y$ is bounded, which will be the case when we apply this result to Huber's estimator but not for Catoni's estimator.

In the rest of the section, we prove concentration inequalities for the estimators featured in Section~\ref{sec:setting_if} using Lemma~\ref{lem:subgaussian} applied to $Y=\frac{X-T(P)}{\|X-T(P)\|}\psi(\|X-T(P)\|)$ and using the bounds on the bias from Section~\ref{sec:bias}. For simplicity, we will not keep track of all the constants and we will give the names $C_1,C_2,C_3$ to numerical constants that do not depend on any of the parameters of the model (in particular they do not depend on $P$ or $\beta$).

\subsection{Huber's estimator}

Let $\beta>0$, and, for all $x\ge 0$, let $\psi_H(x)=x\,\1\{x\le \beta\}+\beta\,\1\{x>\beta\}$.
From Lemma~\ref{lem:ass_exemples}, Assumptions~\ref{ass:prop_psi} hold in this example with $\gamma=1$.
As $\psi_H$ is bounded by $\beta$, Hoeffding's lemma (see \cite[Section 2.3]{concentration}) shows that $\max_{1\le i\le n}\|\psi_H(\|X_i-T_H(P)\|)\|_{\psi_1}\le \beta$. From Lemma~\ref{lem:subgaussian}, we have the following proposition, whose proof is given Section~\ref{sec:proof_huber_corrupted}.

\begin{Proposition}\label{prop:huber_corrupted}
Suppose that Assumption~\ref{ass:adv} and Assumption~\ref{ass:unicity} are verified with $|\mathcal{O}|\le n/32$. Then, there exist some numerical constants $C_1,C_2,C_3>0$ such that if $\beta^2 \ge V_{\psi_H}\max(8,C_1 /n)$, then for any $t \le C_2n$, with probability larger than $1-4e^{-t}-e^{-n/32}$, we have
$$\|T_H(X_1^n)-T_H(P)\|\le 6\sqrt{\frac{V_{\psi_H}}{n}}+8\sqrt{\frac{tv_{\psi_H}}{n}}+\frac{C_3t\beta}{n}+ 8\beta \varepsilon_n,$$

where $\varepsilon_n = |\mathcal{O}|/n$ is the proportion of outliers.
\end{Proposition}
Remark that the condition on $\beta$ can be simplified if needed, using Lemma~\ref{lem:bound_variance_variance}, to $\beta^2 \ge 8Tr(\Sigma)$.
The second step is to choose the value of $\beta$, the choice of $\beta$ will be a tradeoff between the bias term from Lemma~\ref{lem:bias_huber_using_concentration} and the concentration in Proposition~\ref{prop:huber_corrupted}

\begin{Proposition}\label{prop:huber_corrupted_q}
Suppose the same assumptions as in Proposition~\ref{prop:huber_corrupted} and suppose additionally that $\E_P[\|X-\E_P[X]\|^q]<\infty$ for some $q\ge 2$. Then, there exist $C_1, C_2, C_3>0$ numerical constants such that for all $t>0$ with
$$t \le C_1 n\min\left( \frac{\E_P[\|X-\E_P[X]\|^q]}{Tr(\Sigma)^{q/2}}, 1\right),$$
we have with probability larger than $1-4e^{-t}-e^{-n/32}$ that
\begin{multline*}
 \|T_H(X_1^n) - \E_P[X]\|\le 6\sqrt{\frac{Tr(\Sigma)}{n}}+8\sqrt{\frac{t\|\Sigma\|_{op}}{n}}\\
 +C_2\left(\frac{C_3t}{n}+8\varepsilon_n\right)^{1-1/q} \E[\|X-\E[X]\|^q]^{1/q}
 \end{multline*}
where $\beta$ is fixed to the value
$$\beta^q = \frac{2\E[\|X-\E[X]\|^q]}{4\varepsilon_n+C_3t/n}.$$

In particular, if $q>2$ and $t \le C_1 n \frac{\sqrt{\|\Sigma\|_{op}}}{\E[\|X-\E[X]\|^q]^{\frac{1}{q-1}}}$ then we have the following bound on the deviations.
\begin{equation}\label{eq:deviation_hub_simple}
 \|T_H(X_1^n) - \E_P[X]\|\le 6\sqrt{\frac{Tr(\Sigma)}{n}}+9\sqrt{\frac{t\|\Sigma\|_{op}}{n}}+C_2\varepsilon_n^{1-1/q} \E[\|X-\E[X]\|^q]^{1/q}.
\end{equation}
\end{Proposition}
When $\varepsilon_n=0$ the previous proposition guarantees an optimal sub-Gaussian rate.
Notice that $\varepsilon_n$ is multiplied by a quantity that increases with the dimension in general, this bound is not minimax, see~\cite{depersin2019robust} which achieve a sharper bound in the Gaussian setting. We see that the dependency on $\varepsilon_n$ is $O(\varepsilon_n^{1-1/q})$, this type of bound is already present for example in~\cite{DBLP:journals/corr/abs-1903-07870} and the power $1-1/q$ is optimal (see \cite[Lemma 5.4]{minsker2018uniform}) however, the $ \E[\|X-\E[X]\|^q]^{1/q}$ factor in front of it is not optimal, we show in Section~\ref{sec:lower_bound} that this factor is unavoidable for M-estimators. On the other hand, when $\varepsilon_n = 0$, we obtain sub-Gaussian rates of convergence as soon as $t\lesssim n $.

For Equation~\eqref{eq:deviation_hub_simple} to hold, we must have $t \le C_1 n \frac{\sqrt{\|\Sigma\|_{op}}}{\E[\|X-\E[X]\|^q]^{\frac{1}{q-1}}}$. Remark that this is dependent on the dimension, we could have stated a similar deviation bound under the alternative condition  $t \le O\left(n^{\frac{q-2}{2q-2}}\right)$, avoiding the dimension dependence at the price of a worse dependency on $n$. For simplicity we did not state it in the proposition.

Remark that in $\R^d$, with respect to the dimension, $\E[\|X-\E[X]\|^q]^{1/q}$ behaves similarly to $\sqrt{Tr(\Sigma)}$. Indeed, if $\Sigma=I_d$, we have $\sqrt{Tr(\Sigma)}=\sqrt{d}$ and on the other hand, by Jensen's inequality, if we denote by $X^{(i)}$ be the $i^{th}$ coordinate of $X$,
\begin{align*}\E[\|X-\E[X]\|^q]^{1/q} &=\E\left[\left(\sum_{i=1}^d (X^{(i)}-\E[X^{(i)}])^2\right)^{q/2}\right]^{1/q}\\
&\le d^{1/2-1/q}\left(\sum_{i=1}^d \E\left[|X^{(i)}-\E[X^{(i)}]|^q \right]\right)^{1/q}.
\end{align*}
The dependency in the dimension is similarly of order $\sqrt{d}$.

Finally, we present the symmetric case for which there is no need for Lemma~\ref{lem:bias_huber_using_concentration} because we have right-away that $T_H(P)=\E_P[X]$ and this simplifies the computations. In particular we only need a finite first moment and we can directly pick the minimal value of $\beta$ in Proposition~\ref{prop:huber_corrupted} to get,
\begin{Proposition}\label{prop:huber_corrupted_sym}
Suppose the same assumptions as in Proposition~\ref{prop:huber_corrupted} and moreover, suppose that $P$ is symmetric with $\E_P[\|X\|]<\infty$. Then, there exist some numerical constants $C_1, C_2, C_3>0$, such that for any $t \le C_1 n$, with probability larger than $1-4e^{-t}-e^{-n/32}$, we have
\begin{align}
\|T_H(X_1^n)-\E_P[X]\|&\le 6\sqrt{\frac{V_{\psi_H} }{n}}+8\sqrt{\frac{tv_{\psi_H}}{n}}+C_2\E[\|X-\E[X]\|]\left(\frac{t}{n}+ \varepsilon_n\right).
\end{align}
where $\varepsilon_n = |\mathcal{O}|/n$ is the proportion of outliers and $T_H, V_{\psi_H}, v_{\psi_H}$ are computed with $\beta=C_3\E[\|X-\E[X]\|]$.
\end{Proposition}
We see with Proposition~\ref{prop:huber_corrupted_sym} that we can relax the Gaussian inliers assumption made in \cite{depersin2019robust} to inliers that are symmetric or inliers with infinite number of finite moments and still we have a linear dependency in $\varepsilon_n$. We see also that this bound does not need a finite second moment, the law only need to be symmetric and have a finite first moment.

\subsection{Catoni's estimator}

In the case of Catoni's estimator and Polynomial estimator, we will only prove a proposition similar to Proposition \ref{prop:huber_corrupted_q} but we could also state the equivalents of Propositions~\ref{prop:huber_corrupted} and \ref{prop:huber_corrupted_sym} using the same reasoning.

Let $\beta>0$ and, for all $x\ge 0$, let
$\psi_C(x)=
\beta\,\log\left(1+\frac{x}{\beta}+\frac{x^2}{2\beta^2}\right).$
From Lemma~\ref{lem:ass_exemples}, $\psi_C$ satisfies Assumptions~\ref{ass:prop_psi} with $\gamma=4/5$.
Lemma~\ref{lem:subgaussian} in addition to Theorem~\ref{th:tIF_tT_1} can be used to obtain the following proposition.


\begin{Proposition}\label{prop:catoni_corrupted_q}
Suppose Assumption~\ref{ass:adv} and Assumption~\ref{ass:unicity} are verified with $\E_P[\|X\|^2]<\infty$, suppose that there exist some constants $C_O>0$ and $\delta_O \in (0,1)$ such that with probability larger than $1-\delta_O$, we have
$$ \frac{1}{|\mathcal{O}|}\sum_{i\in \mathcal{O}} \psi_C(\|X_i-T_C(P)\|)+\psi_C(\|X'_i-T_C(P)\|) \le \beta C_O.$$
Fix the value of $\beta$ to
$$\beta^3=\frac{5\E[\|X-\E[X]\|^3]}{120C_O\varepsilon_n+Ct\log(n)/n},$$
and suppose $|\mathcal{O}|\le n/(20C_O)$. \\
Then, there exist some numerical constants $C_1,C_2,C_3>0$ such that for any $t\le C_2  \frac{n}{\log(n)^4}\frac{\sqrt{\|\Sigma\|_{op}}}{\E[\|X-\E[X]\|^3]^{1/3}}$, we have with probability larger than $1-4e^{-t}-e^{-n/50}-\delta_O$,
$$ \|T_C(X_1^n) - \E_P[X]\|\le \frac{15}{2}\sqrt{\frac{Tr(\Sigma)}{n}}+11\sqrt{\frac{t\|\Sigma\|_{op}}{n}}+C_2\varepsilon_n^{2/3} \E[\|X-\E[X]\|^3]^{1/3}.$$
\end{Proposition}

Proposition~\ref{prop:catoni_corrupted_q} gives results that are similar to Huber's estimator in Proposition~\ref{prop:huber_corrupted_q} using stronger assumptions on the corruption.
If $\E[\|X\|^3]$ is finite, we can use similar arguments to show a bound of order $O(\varepsilon^{2/3})$ the interested reader could adapt the proof of Proposition~\ref{prop:huber_corrupted_q}. Contrary to Proposition~\ref{prop:huber_corrupted_q}, our method does not allow us to go further than the $3^{rd}$ order because of the bias bound from Lemma~\ref{lem:bias_general}.

The condition on outliers is not very unusual. Indeed, if we suppose that the outliers are i.i.d with law $P_O$ and such that $\E_{P_O}[\psi_C(\|X-T_C(P)\|)]<\infty$, we can use Chebychev inequality to say that for any $C_O>0$
$$\P\left(\forall i \in \mathcal{O}, \quad \psi_C(\|X_i-T_C(P)\|)\le C_O\beta \right)\ge \left(1-\frac{\E_{P_O}[\psi_C(\|X-T_C(P)\|)]}{C_O\beta}\right)^{|\mathcal{O}|} $$
and then, we can choose $C_O$ such that the right-hand-side is strictly positive. Remark that we only suppose a finite moment for $\psi_C(\|X-T_C(P)\|)$, which is a logarithmic moment, this is a rather mild requirement on the outliers. On the other hand, if there is a fixed number of outliers, i.e. $\varepsilon_n = C/n$, we deduce from Proposition~\ref{prop:catoni_corrupted_q} that the speed of convergence does not deteriorate if the outliers are bounded almost surely by $\exp(\sqrt{nV_{\psi_C}})$ in which case $\delta_O=0$. This is also a rather mild requirement on outliers.

%
%

\subsection{Polynomial estimator}
Finally, we look at the polynomial estimator defined for $p,\beta>0$ by
$$\psi_P(x)=\frac{x}{1+(x/\beta)^{1-1/p}}.$$
%



\begin{Proposition}\label{prop:poly_corrupted_2}
Suppose Assumption~\ref{ass:adv} and Assumption~\ref{ass:unicity} are verified with law $P$ that verifies $\E_P[\|X\|^2]<\infty$ and its covariance matrix is denoted $\Sigma$. We suppose that there exist some constants $C_O\ge 1$ and $\delta_O \in (0,1)$ such that with probability $1-\delta_O$, we have
$$ \frac{1}{|\mathcal{O}|}\sum_{i\in \mathcal{O}} \psi_P(\|X_i-T_P(P)\|)+\psi_P(\|X_i'-T_P(P)\|)\le C_O \beta .$$
We fix the value of $\beta$ to $\beta^2=Tr(\Sigma)$.
Then,
there exist some numerical constants $C_1,C_2, C_3>0$ such that for any $t \le C_1 n$, with probability larger than $1-e^{-C_2t}-e^{-n/512}-\delta_O$, we have
$$\|T_P(X_1^n)-\E[X]\|\le 16 \sqrt{65\frac{Tr(\Sigma)t }{n}}+  16\varepsilon_n C_O\sqrt{ Tr(\Sigma)}.
$$
where $\varepsilon_n = |\mathcal{O}|/n$ is the proportion of outliers, supposed positive, and the parameter $p$ is fixed to $p=C_3t$.
\end{Proposition}
The proof of this proposition can be found Section~\ref{sec:proof_poly2} In this proposition we see that we obtain weaker guarantees for the polynomial estimator because the $t$ in the right-hand-side of the bound is multiplied by $Tr(\Sigma)$ as compared to the smaller $\|\Sigma\|_{op}$ factor that we had before, hence this bound is not optimal. Nonetheless, this result is valid with very high probability and this may be surprising to the reader but recall that $p$ is considered a parameter that we tuned using the level $t$ so that in fact we use a function $\psi_P$ that gets very close to a bounded function when $t$ gets large.

%
\section{Lower bound in corruption bias for M-estimators}\label{sec:lower_bound}
In this section we suppose $\psi$ bounded and $\mathcal{H}=\R^d$.
Define $P_\varepsilon=(1-\varepsilon)P+\varepsilon Q$ for some outlier probability $Q$ and some $\varepsilon \in (0,1/2)$. The goal is to estimate the error we would incur if we want to estimate the expectation of $P$ using data from $P_\varepsilon$. Remark that this setting is more restrictive than the adversarial corruption setting we used until now, indeed we can see $P_\varepsilon$ as a corrupted setting in which the adversary chose to modify a random number of randomly chosen outliers.
\begin{Theorem}\label{th:corruption_bias}
Suppose $P=\mathcal{N}(0,\sigma^2 I_d)$. There exist a distribution $Q$ and $\varepsilon_{max}>0$ such that for all $\varepsilon \in (0, \varepsilon_{max})$,
$$\|T(P)-T(P_\varepsilon)\|=\|T(P_\varepsilon)\|\ge \sigma\frac{\varepsilon \sqrt{d-2}}{12}$$
where $P_\varepsilon=(1-\varepsilon)P+\varepsilon Q$.
\end{Theorem}
Theorem~\ref{th:corruption_bias} (proven in Section~\ref{sec:proof_corruption_bias}) gives us a lower bound on the bias due to the corruption. Remark that a similar result already existed for the case of the geometric median, see \cite[Proposition 2.1]{lai2016agnostic}. Our result extends the result of \cite{lai2016agnostic} to more general $\psi$ functions using an alternative proof.

To make the link with the corrupted setting from Assumption~\ref{ass:adv}, remark that if an adversary chose $\lfloor (1-\varepsilon)n\rfloor$ inliers i.i.d from $P$ and $\lceil \varepsilon n\rceil$ outliers i.i.d from $Q$, then the resulting empirical distribution converges to $P_\varepsilon$ and by continuity of $T$ with respect to the weak topology (see Hampel's theorem in~\cite{robuststat}), we have that $T(X_1^n) \to T(P_\varepsilon)$.

Then, remark that we have
\begin{equation}\label{eq:triangle_lower}
\|T(X_1^n)-\E_P[X]\|\ge \|T(P_\varepsilon)-\E[X]\|-\|T(X_1^n)-T(P_\varepsilon)\|.
\end{equation}
In the right-hand-side of \eqref{eq:triangle_lower}, we have the term $\|T(P_\varepsilon)-T(P)\|=\|T(P_\varepsilon)-\E[X]\|$ that stays bounded by below for any value of $\beta$ from Theorem~\ref{th:corruption_bias} and the term $\|T(X_1^n)-T(P_\varepsilon)\|$ that goes to $0$ as $n$ goes to infinity. More precisely, we have the following corollary of Theorem~\ref{th:corruption_bias} consequence of the consistency of the plug-in estimator $T(\widehat P_n)=T(X_1^n)$ for any distribution $P$ (see~\cite{robuststat}).
\begin{Corollary}\label{cor:lower_bound}
Suppose $P=\mathcal{N}(\mu,\sigma^2 I_d)$. There exist a distribution $Q$ and $\varepsilon_{max}>0$ such that for all $\varepsilon \in (0, \varepsilon_{max})$ and if an adversary chose $\lfloor (1-\varepsilon)n\rfloor$ inliers i.i.d from $P$ and $\lceil \varepsilon n\rceil$ outliers i.i.d from $Q$ to form a sample $X_1,\dots,X_n$, we have
$$\P\left(\|T(X_1^n)-\E_P[X]\|\ge  \sigma\frac{\varepsilon \sqrt{d-2}}{24}\right) \xrightarrow[n \to \infty]{} 1.$$
\end{Corollary}
It is possible to quantify the rates in Corollary~\ref{cor:lower_bound} but this is not really necessary as this already proves that we can't hope to achieve a rate that does not depend on $d$ in the corruption error, i.e. we can't attain the minimax rates of convergence which is in this case of order $O(\sigma\varepsilon)$.

\section{M-estimators in practice}\label{sec:algo}
In this section, we give results for $\mathcal{H}=\R^d$ but they could be extended to more general Hilbert spaces provided that one use a sufficiently accurate initialization instead of the coordinate median used here.
\subsection{Algorithm and convergence using iterative re-weighting}
To compute $T(X_1^n)$, we use an iterative re-weighting algorithm. This algorithm is rather well known to compute M-estimators, see \cite[Section 7]{robuststat} and it has already been extensively studied. The principle is to rewrite the definition of $T(X_1^n)$ from equation~\eqref{def:T_hat} as
$$T(X_1^n)\sum_{i=1}^n\frac{\psi\left(\|X_i-T(X_1^n)\|\right)}{\|X_i-T(X_1^n)\|}=\sum_{i=1}^nX_i\frac{\psi\left(\|X_i-T(X_1^n)\|\right)}{\|X_i-T(X_1^n)\|}, $$
then, denote $w_i=\frac{\psi\left(\|X_i-T(X_1^n)\|\right)}{\|X_i-T(X_1^n)\|}$, we get an expression of $T(X_1^n)$ as a weighted sum:
$$T(X_1^n)=\sum_{i=1}^n X_i\frac{w_i}{\sum_{i=1}^n w_i}.$$
The weights $w_i$ depend on $T(X_1^n)$ and the principle of the algorithm is as follows. Initialize $\theta_0$ with the coordinate-wise median and iterate the following
$$w_i^{(m)}= \frac{\psi\left(\|X_i-\theta^{(m)}\|\right)}{\|X_i-\theta^{(m)}\|} \quad \quad \text{,and}\quad \quad  \theta^{(m+1)}=\sum_{i=1}^n X_i\frac{w_i^{(m)}}{\sum_{i=1}^n w_i^{(m)}}.$$
We show that this algorithm allows us to find a minimizer of
$$J_n(\theta)=\frac{1}{n}\sum_{i=1}^n \rho\left(\|X_i-\theta\|\right).$$
Let $r_n,\delta>0$ be such that
$$\P\left(\|T(X_1^n)-\E_P[X]\|\ge r_n \right)\le \delta, $$
for instance, one can use the bound given in Section~\ref{sec:examples_tail}.
We have the following theorem.
\begin{Theorem}\label{th:convergence}
Let $X_1,\dots,X_n$ be in the $\mathcal{I}\cup \mathcal{O}$ setting with $(X_i)_{i \in \mathcal{I}}$ i.i.d with law $P$ whose variance is finite and the covariance matrix is denoted $\Sigma$. Suppose also that $|\mathcal{O}|\le n/8$ and $\beta \ge 2\sqrt{2Tr(\Sigma)}+r_n+\psi^{-1}\left(\sqrt{2V_\psi}\right)$. Then, for all $\N \in \N$, with probability larger than $1-(d+4)e^{-n/8}-\delta$, we have
$$\|\theta^{(m)}-T(X_1^n)\|  \le \frac{1}{(1+\gamma/2)^{m}} \|\theta^{(0)} - T(X_1^n)\| \le \frac{1}{(1+\gamma/2)^{m}} \left(2\sqrt{2Tr(\Sigma)}+r_n \right).$$
Said differently, the iterative reweighting algorithm is such that for any $\varepsilon$, we have $\|\theta^{(m)}-T(X_1^n)\|\le \varepsilon$ after a number of iterations $$m\ge \log\left(\frac{2\sqrt{2Tr(\Sigma)}+r_n }{\varepsilon}\right)/\log(1+\gamma/2).$$
\end{Theorem}
The proof of Theorem~\ref{th:convergence} can be found in Section~\ref{sec:proof_convergence}. To prove this theorem, we use techniques similar to those used to prove the convergence of Weiszfeld’s Method (see~\cite{beck2015weiszfeld}).

We obtain an exponential rate of convergence. Remark that because the objective function is convex, even if the initialization was not as good as the coordinate-wise median or if $\beta$ was not large enough, we would converge nonetheless but with a linear rate of convergence (similar to convergence rates in \cite{beck2015weiszfeld}) until we are close enough to $T(X_1^n)$ for the rate to be exponential.

\subsection{Discussion on the choice of $\beta$}\label{sec:choice_beta}
The choice of $\beta$ is a frequent problem when using Huber estimator. One solution is to use Lepski's method but this is computationally expensive and not always efficient, another (often used) approach is to use a heuristic for $\beta$ based on the median absolute deviation by saying that $\beta$ must be of order $\sigma$ the standard deviation of the inliers however in Section~\ref{sec:examples_tail} we see that this choice is very conservative and $\beta$ would often be too small if it is estimated using the median absolute deviation.

In view of Section~\ref{sec:examples_tail} where we see that depending on the number of finite moments, we may want to choose $\beta$ between $\sqrt{Tr(\Sigma)}$ and $\sqrt{Tr(\Sigma) n}$, we propose to choose $\beta$ in the interval $[0, \mathrm{MAD}\sqrt{n}]$ where $\mathrm{MAD}= \mathrm{Med}\left(\|X_i - \mathrm{GMed}(X_1^n) \| \right)$ (GMED being the geometric median).

We propose the following Heuristic to choose $\beta$:
$$\widehat{\beta} = \argmin_{\beta \in [0, \mathrm{MAD}\sqrt{n}]} \frac{\widehat{V_{\psi_\beta}}}{n} +C_\psi\frac{\mathrm{MAD}^4}{\beta^2} + (0.05\beta)^2 $$
where $\widehat{V_{\psi_\beta}} =\frac{1}{n}\sum_{i=1}^n \psi_\beta(\|X_i-T_{\psi_\beta}(X_1^n)\|)^2$ and $C_\psi$ is a constant that depend on $\psi$ as described by the bounds on the bias (Section~\ref{sec:bias}, $C_{\psi_H} = 1$, $C_{\psi_C} = 5/32$ and $C_{\psi_P}=1/16$). This is a bias-variance trade-off, the first term converges to the asymptotic variance of $T_{\psi_\beta}(X_1^n)$, the second term is a bound on the squared bias and the third term is a bound on the corruption bias if we suppose that $\varepsilon_n \le 0.05$ (In the Robust literature, it is often said that there is less than $5$ or $10$ percent of outliers). 

Remark that the objective function may not be convex and hence there can be local minima, we restrict the search space $[0, \mathrm{MAD}\sqrt{n}]$ in order to be able to choose $\beta$ efficiently using a grid-search.

\subsection{Illustrations}\label{sec:illustration}
To illustrate the behavior of M-estimators in heavy-tailed setting we consider a multivariate Pareto law for which the coordinates are drawn, independently of each other, from a Pareto distribution and a multivariate student distribution. All of the dataset present a finite variance but infinite third moment.
\begin{description}
\item[Dataset 1]: Coordinates drawn from a Pareto with shape parameter $\alpha=2.1$ and scale parameter $1$, with two corrupted samples situated in $300 \textbf{1}_d$ where $\textbf{1}_d$ is the vector with all coordinates equal to $1$.
\item[Dataset 2]: Coordinates drawn from a Pareto with shape parameter $\alpha=3$ and scale parameter $1$. Uncorrupted dataset.
\item[Dataset 3]: Data i.i.d from a mixture $0.4 \mathcal{T}(0, 2.1) + 0.6 \mathcal{T}(2\textbf{1}_d,2.1)$, where $\mathcal{T}(\mu, \nu)$ is the multivariate student distribution with mean $\mu$ and degree of freedom $\nu$. In this dataset, there is also two corrupted samples situated in $300 \textbf{1}_d$.
\item[Dataset 4]: Data i.i.d from the mixture $0.4 \mathcal{T}(0, 3) + 0.6 \mathcal{T}(2\textbf{1}_d,3)$. Uncorrupted dataset.
\end{description}
In these three datasets, we consider $n=1000$ samples and the dimension is $d=100$.

\begin{figure}[h]
  \begin{center}
    \subfloat[Dataset 1]{
\includegraphics[scale=0.4]{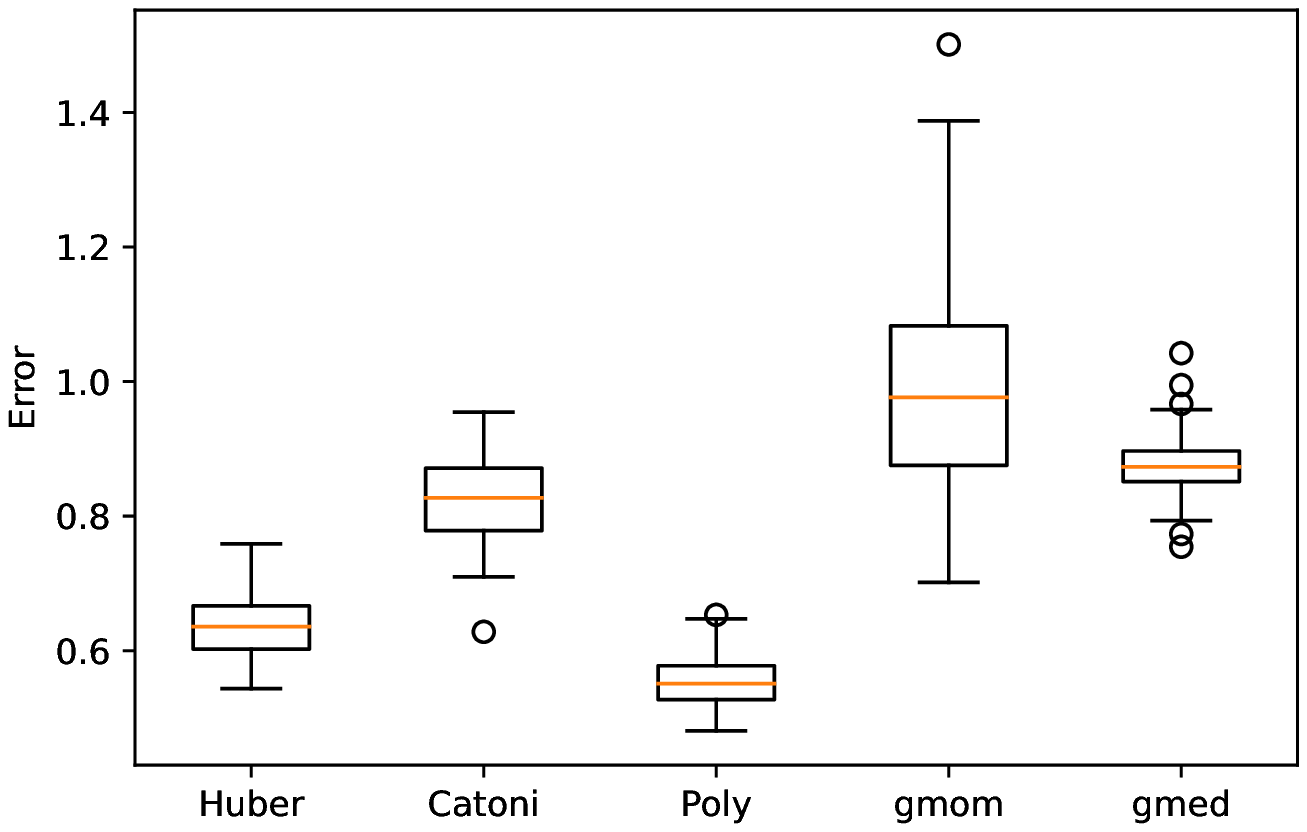}}
\subfloat[Dataset 2]{
\includegraphics[scale=0.4]{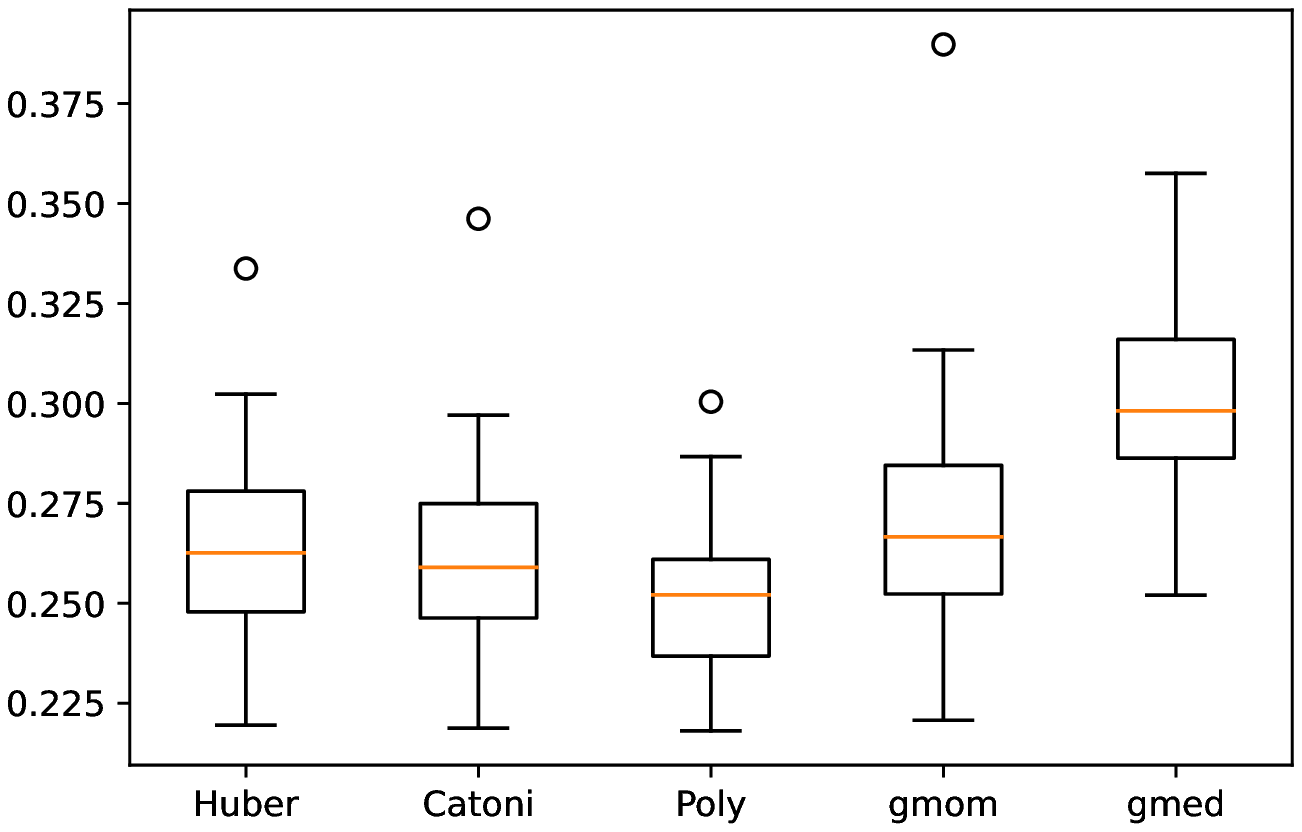}}\\
\subfloat[Dataset 3]{
\includegraphics[scale=0.4]{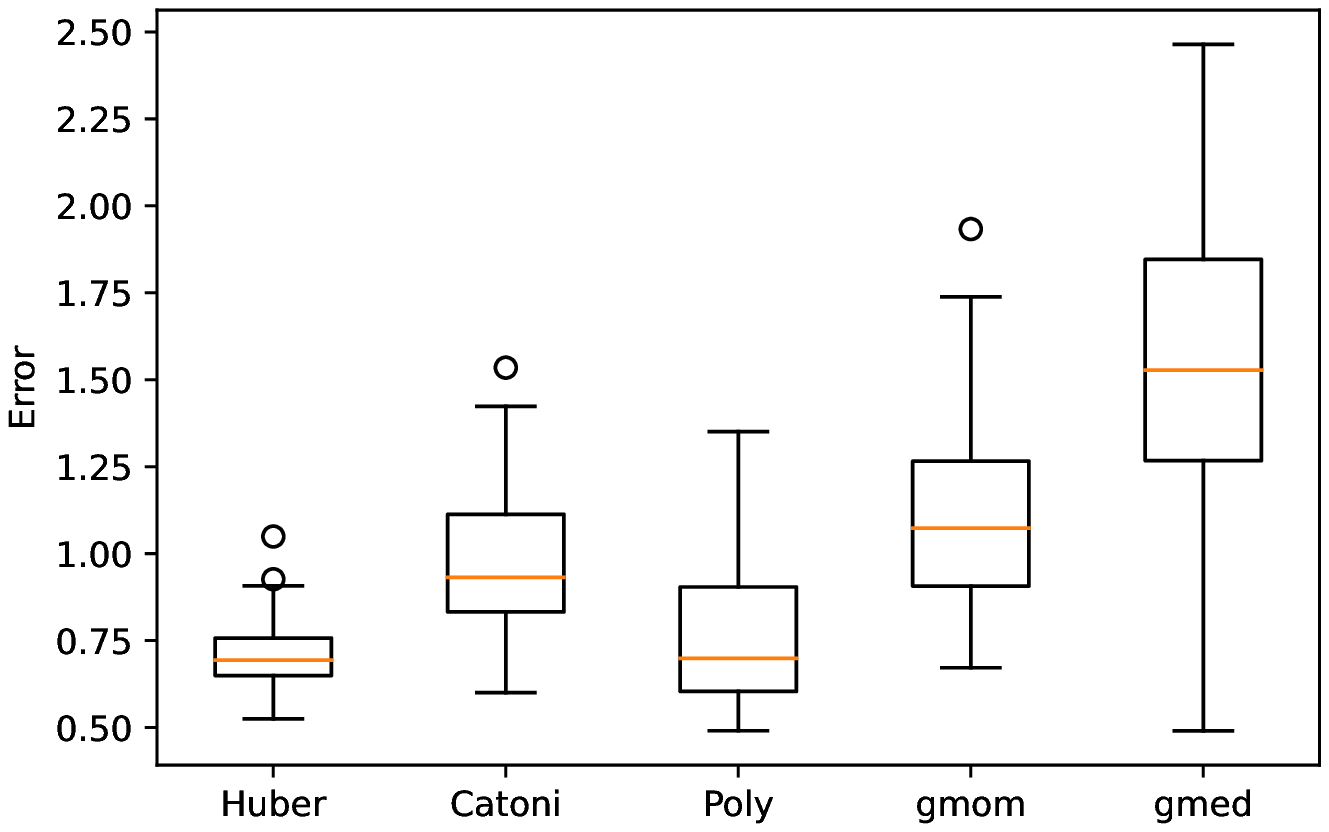}}
\subfloat[Dataset 4]{
\includegraphics[scale=0.4]{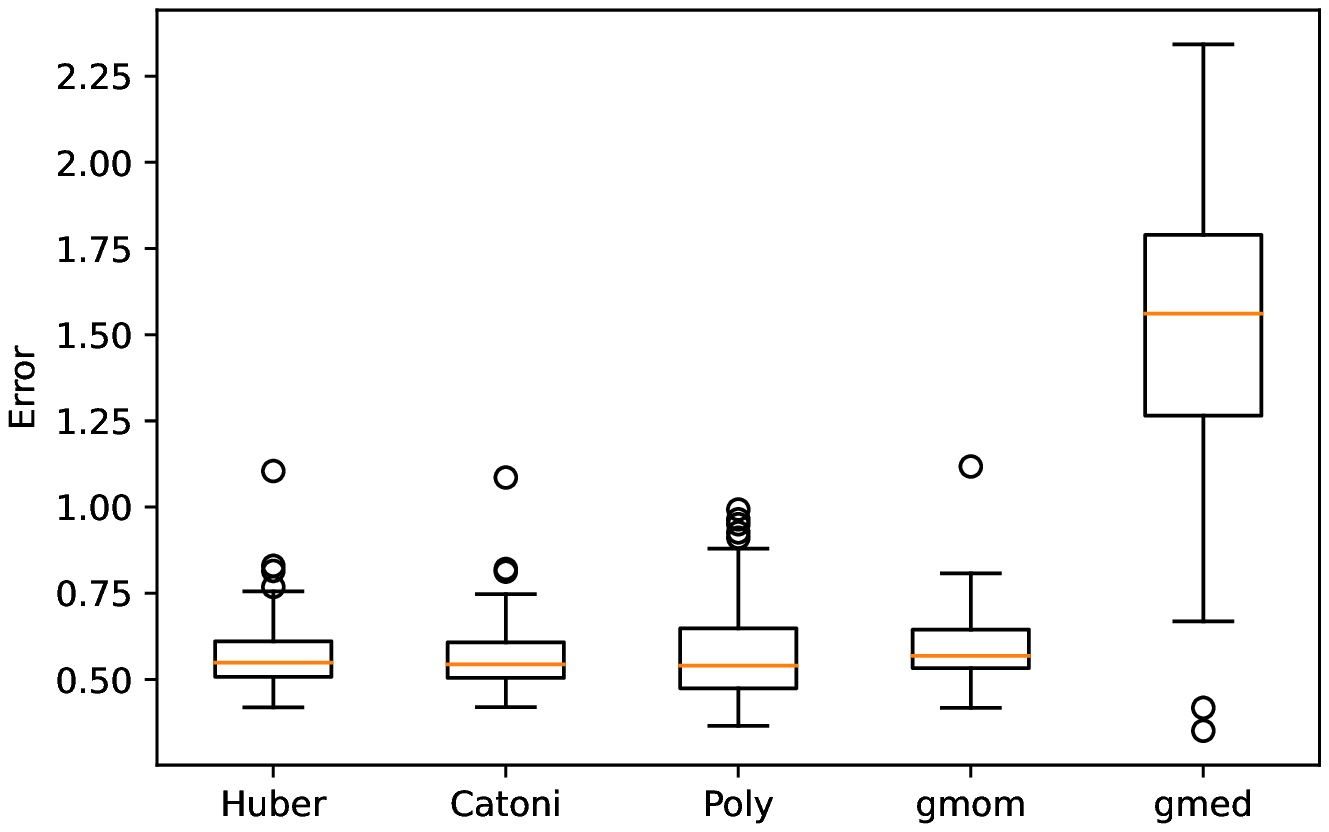}}
\caption{Boxplots of the error $\|\widehat{\mu} - \E[X]\|$ of various estimators $\widehat{\mu}$ using $100$ independent copies of the three dataset considered.\label{fig:illustration}}
\end{center}
\end{figure}

We consider 5 estimators : Huber and Catoni estimator with $\beta$ chosen using Section~\ref{sec:choice_beta}, the polynomial estimator with $p=5$ and $\beta$ chosen using Section~\ref{sec:choice_beta}, the geometric median (gmed) and the geometric median of means described in \cite{MR3378468} with $k=9$ blocks. The result is represented in Figure~\ref{fig:illustration}

In Figure~\ref{fig:illustration}, M-estimators are outperforming the geometric median by quite a lot because the geometric median is very biased when estimating the mean of an asymmetric distribution. The geometric median of means is closer in performance to M-estimators but it is not as good, maybe because there is no adaptivity the number of blocks have been fixed to $9$..

Remark that the multivariate robust mean estimators described in \cite{diakonikolas2020outlier,depersin2019robust,hopkins2020mean,cherapanamjeri2019fast} are not presented here because they are either too computationally intensive or too hard to implement for the purpose of comparison. Remark also that the choice of $\beta$ from Section~\ref{sec:choice_beta} is heuristic. In a learning setting, one may prefer to use cross-validation to tune $\beta$ using directly the learning criterion.

\begin{appendix}
\section{Proof of Theorems}
\subsection{Proof of Theorem \ref{th:tIF_tT_1}}
\label{sec:proof_tIF_tT}
\begin{Hypopf}
\begin{enumerate}[label=(\roman{*}), ref=Hypothesis (\roman{*}), leftmargin=*]
\item $V_\theta=\E_P[\psi(\|X'-\theta\|)^2]\le \psi(\beta/2)^2/2<\infty$\label{hypopf:V}
\item $|\mathcal{O}|\le n\gamma/8$\label{hypopf:outliers}
\end{enumerate}
\end{Hypopf}

For all $n\in\N^*$, $\lambda\in \R, \theta \in \mathcal{H}$ and $\theta \neq T(X_1^n)$, let
$$f_{n}(\lambda; \theta)=\frac{1}{n}\sum_{i=1}^n \frac{\langle X_i-\theta-\lambda u_n(\theta),u_n(\theta) \rangle }{\|X_i-\theta-\lambda u_n(\theta)\|} \psi(\|X_i-\theta-\lambda u_n(\theta)\|),$$
where $u_n(\theta)= \frac{\theta - T(X_1^n)}{\|\theta - T(X_1^n)\|}$.
\begin{npf}
\item For any $\lambda>0$, $\P(\| \theta-T(X_1^n)\| \ge \lambda)\le \P( f_{n}(\lambda; \theta)\ge 0).\label{item:goal}
$
\\
\onepf{
For all $\theta \in \mathcal{H}$, let $J_n(\theta)=\frac{1}{n}\sum_{i=1}^n \rho(\|X_i-\theta\|)$  from Lemma~\ref{lem:convexity} applied to the empirical density $\widehat{P}_n = \frac{1}{n}\sum_{i=1}^n \delta_{X_i}$, we have for any $u \in S$ and $\theta \in \mathcal{H}$,
\[
u^T Hess(J_n)(\theta)u \ge \frac{1}{n}\sum_{i=1}^n \psi'\left(\|X_i-\theta\|\right).
\]
In particular, $f_{n}(\lambda;\theta)=-\langle \nabla J_n(\theta - \lambda u_n(\theta)), u_n(\theta)\rangle$ and if we take the derivative of $f_{n}$ with respect to $\lambda$, we have the following equation
\begin{equation}\label{eq:deriv_neg}
f_{n}'(\lambda;\theta)=-u_n(\theta)^T Hess(J_n)(\theta-\lambda u_n(\theta))u_n(\theta)\le-\frac{1}{n}\sum_{i=1}^n\psi'(\|X_i-\theta-\lambda u_n(\theta)\|).
\end{equation}
where $f_{n}'(\lambda;\theta)$ is the derivative of $f_{n}(\lambda;\theta)$ with respect to $\lambda$.
Then, because $\psi'$ is non-negative (Assumption~\ref{ass:3}), the function $f_{n}(\,\cdot\,; \theta)$ is non-increasing. Hence, for all $n\in\N^*$, $\lambda\in \R$
$$ \| \theta-T(X_1^n)\| \ge \lambda \Rightarrow  f_{n}( \| \theta -T(X_1^n)\| ; \theta)=0 \le f_{n}(\lambda;\theta).$$
And then,
\begin{align}\label{eq:oracle1}
\P(\| \theta-T(X_1^n)\| \ge \lambda)&\le \P( f_{n}(\lambda; \theta)\ge 0).\end{align}}
\item For all $\lambda>0$, \label{item:taylor}
\begin{equation*}
f_{n}(\lambda; \theta)\le f_{n}(0)-\lambda\inf_{t \in [0,\lambda]}\left|f_{n}'(t;\theta)\right|. \end{equation*}\\
\onepf{
We apply Taylor's inequality to the function $\lambda \mapsto f_{n}(\lambda; \theta)$.
As $\lambda \mapsto f_{n}(\lambda; \theta)$ is non-increasing (because its derivative is non-positive, see Equation~\eqref{eq:deriv_neg}), we get the result.
}
\item Let $m = \gamma \P\left( \|X_i'-\theta\|\le \beta/2\right)$. With probability  than $1-e^{-n m^2/8}$, we have \label{item:stp3}
$$\inf_{ t \in [0,\lambda]}|f_{n}'(\lambda; \theta)|\ge  \frac{3m}{4}-\frac{|\mathcal{O}|}{n}.$$
\onepf{
From equation~\eqref{eq:deriv_neg},
\begin{align}\label{eq:oracle2}
|f_{n}'(t;\theta))|\ge& \frac{1}{n}\sum_{i=1}^n \psi'(\|X_i-\theta-t u_n(\theta)\|)\nonumber\\
\ge & \frac{1}{n}\sum_{i=1}^n \psi'(\|X_i'-\theta-t u_n(\theta)\|) \nonumber \\
&+ \frac{1}{n}\sum_{i \in \mathcal{O}} \left( \psi'(\|X_i-\theta-t u_n(\theta)\|)-\psi'(\|X_i'-\theta-t u_n(\theta)\|) \right).
\end{align}
Hence,  because $\psi' \in [0,1]$, we have
\begin{align*}
\inf_{t \in [0,\lambda]}|f_{n}'(t;\theta))| &\ge \inf_{t \in [0,\lambda]} \frac{1}{n}\sum_{i=1}^n \psi'(\|X_i'-\theta-t u_n(\theta)\|) - \frac{|\mathcal{O}|}{n}\\
&\ge \inf_{\substack{ t \in [0,\lambda]}} \frac{\gamma}{n}\sum_{i=1}^n \1\{\|X_i'-\theta-t u_n(\theta)\|\le \beta\} - \frac{|\mathcal{O}|}{n}\\
&\ge \frac{\gamma}{n}\sum_{i=1}^n \1\{\|X_i'-\theta\|\le \beta/2\} - \frac{|\mathcal{O}|}{n}
\end{align*}
The last line is a consequence of $\lambda \in (0,\beta/2)$. Then, this is a sum of bounded random variables and by Hoeffding's inequality, with probability larger than $1-e^{-2n\varepsilon^2}$, we have
$$\sup_{ t \in [0,\lambda]}-|f_{n}'(\lambda; \theta)|\le -\gamma\E\left[ \1\{\|X_i'-\theta\|\le \beta/2\}\right]+\varepsilon+\frac{|\mathcal{O}|}{n}.$$
take $\varepsilon = \gamma \P\left( \|X_i'-\theta\|\le \beta/2\right)/4$ to conclude
}
\item For $\lambda\in (0, \beta/2)$, $\P(\| \theta-T(X_1^n)\| \ge \lambda)\le  t_{\IF}\left(\lambda\gamma/4;\theta, X_1^n\right)+e^{-n\gamma^2/32}.$\\
\onepf{
For any $\lambda>0$, we have
\begin{align}\label{eq:result_concentration}
\P(\| \theta-T(X_1^n)\| \ge \lambda)&\stackrel{(a)}{\le} \P( f_{n}(\lambda;\theta)\ge 0)\nonumber\\
&\stackrel{(b)}{\le} 1-\P\left(f_{n}(0;\theta)-\lambda\inf_{t \in [0,\lambda]}\left|f_{n}'(t;\theta)\right|\le 0\right)\nonumber\\
&\stackrel{(c)}{\le} 1-\P\left( f_{n}(0;\theta)\le \lambda\left(\frac{3}{4}m - \frac{|\mathcal{O}|}{n} \right) \right)+e^{-n m^2/8}\nonumber\\
&= t_{\IF}\left(\lambda \left(\frac{3}{4}m - \frac{|\mathcal{O}|}{n} \right);\theta,X_1^n\right)+e^{-n m^2/8},
\end{align}
where $(a)$ follows from \ref{item:goal}, $(b)$ follows from \ref{item:taylor} and $(c)$ follows from \ref{item:stp3}.
Finally, we bound $m$ because $\psi$ is increasing on $[0,\beta]$, for any $\lambda < \beta/2$,
\begin{align*}
 \P\left(\|X_i'-\theta\|\le \beta/2 \right)= \P\left(\psi(\|X_i'-\theta\|)\le \psi(\beta/2) \right).
\end{align*}
As $V_\theta=\E[\psi(\|X-\theta\|)^2]<\infty$, by Markov's inequality, it follows that
\begin{align}\label{eq:markov_upper}
m=\gamma \P\left(\|X_i'-\theta\|\le \beta/2 \right)&\ge  \gamma \left(1-\frac{V_\theta}{\psi(\beta/2)^2} \right).
\end{align}
From~\ref{hypopf:V}, $V_\theta\le \psi(\beta/2)^2/2$, we get $m\ge \gamma/2$.
Plugging this bound in equation \eqref{eq:result_concentration}, we get
\begin{align}
\P(\| \theta-T(X_1^n)\| \ge \lambda)&\le  t_{\IF}\left(\lambda \left(3\gamma/8 - \frac{|\mathcal{O}|}{n} \right);\theta, X_1^n\right)+e^{-n\gamma^2/32}.
\end{align}
we conclude using that $|\mathcal{O}|\ge n\gamma/8$ (from \ref{hypopf:outliers}).}
\end{npf}
\subsection{Proof of Theorem~\ref{th:cornerstone}}\label{sec:proof_cornerstone}
\begin{npf}
\item We have
\begin{multline*}\left\|Z_\beta(\E[X])-Z_\beta(T(P))\right\|\\\ge \gamma\|\E[X]-T(P)\|  \P\big(\left\|X-T(P)\right\|\le \beta-\|\E[X]-T(P)\|\big).
\end{multline*}\label{item:cornerstone:1}\\
\onepf{
The function $\theta \mapsto Z_\beta(\theta)$ is differentiable and by the mean value theorem, we have
\begin{multline}\label{eq:lower_taylor}
\left\|Z_\beta(\E[X])-Z_\beta(T(P))\right\| \\\ge \|\E[X]-T(P)\|\inf_{t\in [0,1]}\left\|Jac( Z_\beta)(t\E[X]+(1-t)T(P))\right\|_{op}
\end{multline}
Where $Jac$ denotes the Jacobi operator. From Lemma~\ref{lem:convexity} and Assumption~\ref{ass:3}, we get
$$\left\|Jac(Z_\beta)(\theta)\right\|_{op}\ge \E\left[\psi_1'\left(\left\|\frac{X-\theta}{\beta}\right\|\right) \right]\ge \gamma \P\left(\left\|X-\theta\right\|\le \beta\right) .$$
Hence, for all $t\in [0,1]$,
\begin{align*}
\left\|Jac(Z_\beta)(t\E[X]+(1-t)T(P))\right\|_{op}&\ge  \gamma \P\left(\left\|X-t\E[X]-(1-t)T(P)\right\|\le \beta\right) \\
&\ge \gamma \P\left(\left\|X-T(P)\right\|\le \beta-\|\E[X]-T(P)\|\right)
\end{align*}
inject this in Equation~\eqref{eq:lower_taylor} to get the result.}
\item $\P\left(\left\|X-T(P)\right\|\le \beta-\|\E[X]-T(P)\|\right)\ge \frac{1}{2}.$\label{item:cornerstone:2}\\
\onepf{
Use the following lemma proven in Section~\ref{sec:proof_rough}.
\begin{Lemma}\label{lem:rough}
If $\rho(1/3)\ge \E[\rho(\|X-\E[X]\|/\beta)]$, then
$\|\E[X]-T(P)\|\le \frac{\beta}{3} $.
\end{Lemma}
We get,
\begin{align*}
\P\left(\left\|X-T(P)\right\|\le \beta-\|\E[X]-T(P)\|\right)&\ge  \P\left(\left\|X-T(P)\right\|\le 2\beta/3\right)\\
&= \P\left(\rho\left(\frac{\left\|X-T(P)\right\|}{\beta }\right)\le \rho(2/3)\right)\\
&\ge  \P\left(\rho\left(\frac{\left\|X-T(P)\right\|}{\beta }\right)\le 2\rho(1/3)\right)\\
\end{align*}
because $\rho$ is increasing and super-additive on $\R_+$ ($\rho$ is increasing because $\psi(0)=0$ and $\psi$ is non-decreasing because $\psi'\ge 0$, hence $\psi=\rho'\ge 0$). Hence, by Markov's inequality and using the hypothesis,
\begin{multline*}
\P\left(\left\|X-T(P)\right\|\le \beta-\|\E[X]-T(P)\|\right)\\
\ge 1- \frac{1}{2\rho(1/3)}\E\left[\rho\left(\frac{\left\|X-T(P)\right\|}{\beta}\right)\right]\ge \frac{1}{2}
\end{multline*}
}
\item $ \|\E[X]-T(P)\|\le \frac{2}{\gamma}\left\|Z_\beta(\E[X])\right\|.$\\
\onepf{
Use \ref{item:cornerstone:1} and \ref{item:cornerstone:2}}
\end{npf}

\subsection{Proof of Theorem~\ref{th:corruption_bias}}\label{sec:proof_corruption_bias}
$P$ is symmetric around $0$, hence $T(P)=0$. By definition of the influence function we have
$$\IF(x)=\lim_{\varepsilon \to 0}\frac{T(P_\varepsilon)-T(P)}{\varepsilon} $$
where $P_\varepsilon=(1-\varepsilon)P+\varepsilon \delta_x$
and we have, from~\cite{HampelEtal86} that
$\IF(x)=M^{-1}\frac{x}{\|x\|}\psi(\|x\|),$
where
$$
M=-\E\left[\left(\frac{I_d}{\|X\|}-\frac{XX^T}{\|X\|^3} \right)\psi(\|X\|)+\frac{XX^T}{\|X\|^2}\psi'(\|X\|)  \right] .
$$
Hence, taking $\|x\|$ such that $\psi(\|x\|)\ge \|\psi\|_\infty/2$ (this fixes the norm of $x$) and $\|M^{-1}x\|=\|M^{-1}\|_{op}\|x\|$ (this fixes the direction of $x$), we have
\begin{equation}\label{eq:bound_if_proof}
\|\IF(x)\|\ge \frac{\|M^{-1}\|_{op}}{2}\|\psi\|_\infty\ge \frac{1}{2\|M\|_{op}}\|\psi\|_\infty
\end{equation}
where $\|M\|_{op}$ is the operator norm of $M$. Let us control this operator norm. We have for all $u \in S$ where $S$ is the sphere in $\R^d$,
$$u^TMu = -\E\left[\left(\frac{1}{\|X\|}-\frac{\langle X, u\rangle ^2 }{\|X\|^3} \right)\psi(\|X\|)+\frac{\langle X, u\rangle ^2}{\|X\|^2}\psi'(\|X\|)  \right] $$
Hence,
$$|u^TMu| \le  3\E\left[\frac{\psi(\|X\|)}{\|X\|} \right]\le 3 \|\psi\|_\infty \E\left[\frac{1}{\|X\|} \right] \le 3 \frac{\|\psi\|_\infty}{\sigma} \E\left[\frac{\sigma^2 }{\|X\|^2} \right]^{1/2} $$
Then, use that $X \sim \mathcal{N}(0,\sigma^2 I_d)$ to have that $\sigma^2/\|X\|^2$ have a law inverse-$\chi^2$ of parameter $d$. Hence,
$$|u^TMu| \le  3 \frac{\|\psi\|_\infty }{\sigma \sqrt{d-2}}.$$
Inject this in equation~\eqref{eq:bound_if_proof} to conclude that $\|\IF(x)\|\ge \sqrt{d-2}/6.$
 Now, use that $\varepsilon \mapsto T(P_\varepsilon)$ is continuous (it is in fact Lipshitz continuous because the Influence function is bounded) to conclude.
\subsection{Proof of Theorem~\ref{th:convergence}}\label{sec:proof_convergence}
This proof mimic the proof of convergence for EM algorithm or for the algorithm used to compute the geometric median.

Let
$$J_n(\theta) = \frac{1}{n}\sum_{i=1}^n \rho_1\left(\frac{\|X_i-\theta\|}{\beta} \right),$$ we are searching for the argmin of $J_n$.

First, we show that the initialization is not too far away from the optimum.
\begin{Lemma}\label{lem:initialization}
If Assumption~\ref{ass:adv} is verified, $\Sigma$ is the covariance matrix of $P$ and $|\mathcal{O}|\le n/8$.
Then, with probability larger than $1-de^{-n/8}-\delta$,
$$\|\theta_0 - T(X_1^n)\|\le 2\sqrt{2Tr(\Sigma)}+r_n$$
\end{Lemma}
The proof is provided in Section~\ref{sec:proof_lem_initialization}. Then, we note that $J_n$ is strongly convex with high probability, using Lemma~\ref{lem:convexity},
for all $u \in S$ and all $\theta \in \R^d$,
\begin{equation}\label{eq:strongly_convex}
u^T Hess(J_n)(\theta)u \ge \frac{1}{\beta^2 n}\sum_{i=1}^n \psi_1'\left(\frac{\|X_i-\theta\|}{\beta}\right).
\end{equation}
This allows us to show that the sequence of iterates is equivalent to minimizing a convex majorant of $J$.
\begin{Lemma}\label{lem:surrogate}
Assume that $\rho$ is convex, that $x \mapsto \psi(x)/x$ is bounded and non-increasing, define
$$\textbf{U}_\kappa(\theta)=\sum_{i=1}^n\left( \frac{w_i(\kappa)}{2}\left(\frac{\|X_i-\theta\|}{\beta} \right)^2+ \rho_1(r_i(\kappa))-\frac{1}{2}r_i(\kappa)\psi(r_i(\kappa))\right). $$
where $r_i(\kappa)= \|X_i-\kappa\|/\beta$ and $w_i(\kappa)=\psi_1(r_i(\kappa))/r_i(\kappa)$. We have that
\begin{itemize}
\item  $\textbf{U}_{\theta^{(m)}}$ is minimized at $\theta^{(m+1)}$
\item for all $\kappa \in \R^d$, $\textbf{U}_\kappa\ge J_n$
\item The approximation error $h=J_n-U_\kappa$ is differentiable and $\nabla h$ is $L$-Lipshitz continuous with $L=\frac{1}{n\beta^2}\sum_{i=1}^n \frac{\psi(r_i)}{r_i}$
\item $h(\kappa)=0$ and $\nabla h(\kappa)=0$.
\end{itemize}

\end{Lemma}
The proof is provided in Section~\ref{sec:proof_lem_surrogate}.

From Lemma~\ref{lem:surrogate}, we get that the direction from $\theta^{(m)}$ to $\theta^{(m+1)}$ is a proper descent direction through the following lemma proved in Section~\ref{sec:proof_lem_descent}.
\begin{Lemma}\label{lem:descent}
For all $\theta \in \R^d$ and for all $m \in \N$, we have
$$J(\theta^{(m+1)}) \le  J(\theta)+\frac{1}{2\beta^2n}\sum_{j=1}^n w_j(\theta^{(m)}) \left(\|\theta^{(m)} - \theta\| - \|\theta^{(m+1)}-\theta\|\right) $$
\end{Lemma}

Use Lemma~\ref{lem:descent} with $\theta=T(X_1^n)$.
\begin{multline*}J(\theta^{(m+1)}) \le  J(T(X_1^n))\\+\frac{1}{2\beta^2n}\sum_{j=1}^n w_j(\theta^{(m)}) \left(\|\theta^{(m)} - T(X_1^n)\| - \|\theta^{(m+1)}-T(X_1^n)\|\right) 
\end{multline*}
Hence, using the fact that $T(X_1^n)$ is a minimizer of $J$, we get the result that
\begin{equation}\label{eq:decreasing_sequence}
\|\theta^{(m+1)}-T(X_1^n)\| \le \|\theta^{(m)} - T(X_1^n)\|
\end{equation}

This allows us to restrict ourselves to a bounded domain. From Lemma~\ref{lem:initialization}, with probability larger than $1-de^{-n/8}-\delta$, we have that $\theta_0$ is in
$$\Theta = \{\theta \in \R^d:\quad \|\theta - T(X_1^n)\|\le 2\sqrt{2 Tr(\Sigma)}+r_n\}.$$
and then, equation~\eqref{eq:decreasing_sequence} assures us that we stay in $\Theta$ for the other iterations: $\forall m \in \N, \theta^{(m)} \in \Theta$. The following Lemma is proven in Section~\ref{sec:proof_convex2}.

\begin{Lemma}\label{lem:convex2}
Let $\theta \in \Theta$, if $X_1,\dots,X_n$ are corrupted by an adversary with $X_1',\dots,X_n'$ i.i.d with law $P$ whose variance is finite and the covariance matrix is denoted $\Sigma$. Moreover, we suppose $|\mathcal{O}|\le n/8$,
 then if $\beta\ge 2\sqrt{2Tr(\Sigma)}+r_n+\psi_1^{-1}\left(\sqrt{2V_\psi}\right)$, we have with probability greater than $1-e^{-n/32}$,
$$\frac{1}{\beta^2n}\sum_{i=1}^n \psi_1'\left(\frac{\|X_i-\theta\|}{\beta}\right)\ge \frac{\gamma}{4\beta^2}$$
\end{Lemma}
This allows us to quantify the speed at which the sequence $\|\theta^{(m)}-T(X_1^n)\|$ decreases. Indeed, from Lemma~\ref{lem:descent} for $\theta=T(X_1^n)$, we have
\begin{multline*}
J(\theta^{(m+1)}) \le  J(T(X_1^n))\\+\frac{1}{2\beta^2n}\sum_{j=1}^n w_j(\theta^{(m)}) \left(\|\theta^{(m)} - T(X_1^n)\| - \|\theta^{(m+1)}-T(X_1^n)\|\right) 
 \end{multline*}
and by Lemma~\ref{lem:convex2} and Lemma~\ref{lem:initialization} and the convexity equation~\eqref{eq:strongly_convex}, we have by Taylor's theorem that with probability larger than $1-(d+4)e^{-n/8}-\delta$ (because $e^{-n/32}\le 4e^{-n/8}$),
$$J(\theta^{(m+1)})-J(T(X_1^n)) \ge \frac{\gamma}{4\beta^2} \|\theta^{(m+1)}-T(X_1^n)\|.$$ Hence,
\begin{multline*}
\frac{\gamma}{4\beta^2} \|\theta^{(m+1)}-T(X_1^n)\| \le \\\frac{1}{2\beta^2n}\sum_{j=1}^n w_j(\theta^{(m)}) \left(\|\theta^{(m)} - T(X_1^n)\| - \|\theta^{(m+1)}-T(X_1^n)\|\right).
\end{multline*}
Solve this equation for $\|\theta^{(m+1)}-T(X_1^n)\| $ to obtain,
\begin{equation}\label{eq:consistance1}
\|\theta^{(m+1)}-T(X_1^n)\|  \le \frac{\frac{1}{2n\beta^2}\sum_{j=1}^n w_j(\theta^{(m)})}{\frac{1}{2n\beta^2}\sum_{j=1}^n w_j(\theta^{(m)})+\frac{\gamma}{4\beta^2}} \|\theta^{(m)} - T(X_1^n)\|
\end{equation}
then use that $\sum_{j=1}^n w_j(\theta^{(m)}) \le n$ and the fact that the right-hand-side of Equation~\eqref{eq:consistance1} is increasing in $\sum_{j=1}^n w_j(\theta^{(m)})$ to get for all $m\in \N$,
\begin{equation}\label{eq:consistance2}
\|\theta^{(m+1)}-T(X_1^n)\|  \le \frac{1/2}{1/2+\frac{\gamma}{4}} \|\theta^{(m)} - T(X_1^n)\|=\frac{1}{1+\gamma/2} \|\theta^{(m)} - T(X_1^n)\|
\end{equation}
Then, this implies directly that with probability larger than $1-(d+4)e^{-n/8}-\delta$, we have
$$\|\theta^{(m)}-T(X_1^n)\|  \le \frac{1}{(1+\gamma/2)^{m}} \|\theta^{(0)} - T(X_1^n)\| \le \frac{1}{(1+\gamma/2)^{m}} \left(2\sqrt{2Tr(\Sigma)}+r_n \right)$$
\section{Proof of Propositions}
\subsection{Proof of Proposition~\ref{prop:huber_corrupted}}\label{sec:proof_huber_corrupted}
\begin{Hypopf}
\begin{enumerate}[label=(\roman{*}), ref=Hypothesis (\roman{*}), leftmargin=*]
\item $|\mathcal{O}|\le n/32$ \label{item:hyp_prop1:1}
\item $\beta^2 \ge V_{\psi_H}\max(8, C_1/n)$\label{item:hyp_prop1:2}
\end{enumerate}
\end{Hypopf}

From Theorem~\ref{th:tIF_tT_1} with $\gamma=1$, we have that if $V_{\psi_H}=\E_P[\psi_H(\|X-T_H(P)\|)^2]\le \psi(\beta/2)^2/2=\beta^2/8$, then for all $\lambda \in (0,\beta/2)$ and for all $\theta \in H$,
\begin{equation}\label{eq:tif_tt_h}
t_T(\lambda; T_H(P), X_1^n)\le t_{\IF}\left(\lambda/4; T_H(P), X_1^n\right)+e^{-n/32}.
\end{equation}
let us find an upper bound of $t_{\IF}\left(\lambda/4; T_H(P), X_1^n\right)$.
\begin{npf}
\item $t_{\IF}\left(\lambda_t/4; T_H(P),X_1^n\right)\le  4e^{-t}.$\\
\onepf{We have,
\begin{multline}\label{eq:tifh1}
\P\left(\left\|\frac{1}{n}\sum_{i=1}^n \frac{X_i-T(P)}{\|X_i-T_H(P)\|}\psi_H\left(\|X_i-T_H(P)\|\right) \right\| \ge \lambda /4\right)\\
\le \P\left(\left\|\frac{1}{n}\sum_{i=1}^n \frac{X_i'-T(P)}{\|X_i'-T_H(P)\|}\psi_H\left(\|X_i'-T_H(P)\|\right) \right\| \ge \lambda /4 - 2\frac{|\mathcal{O}|}{n} \beta \right)
\end{multline}
because $\psi_H$ is bounded by $\beta$. Then, we use Lemma~\ref{lem:subgaussian} to bound the sum of i.i.d. bounded random variables. Having $\psi_H$ bounded, we see that by Hoeffding's lemma (see \cite[Section 2.3]{concentration}), $\max_{1\le i\le n}\|\psi_H(\|X_i-T_H(P)\|)\|_{\psi_1}\le \beta$. Then, by Lemma~\ref{lem:subgaussian} for all $t>0$, with probability larger than $1-4e^{-t}$,
$$\left\|\frac{1}{n}\sum_{i=1}^n \frac{X_i'-T(P)}{\|X_i'-T_H(P)\|}\psi_H\left(\|X_i'-T_H(P)\|\right) \right\|\le  \frac{3}{2}\sqrt{\frac{V_{\psi_H}}{n}}+2\sqrt{\frac{tv_{\psi_H}}{n}}+\frac{Ct\beta}{n} ,$$
inject this in Equation~\eqref{eq:tifh1} with $\lambda_t = 4 \left(\frac{3}{2}\sqrt{\frac{V_{\psi_H}}{n}}+2\sqrt{\frac{tv_{\psi_H}}{n}}+\frac{Ct\beta}{n} \right)+ 8\beta \frac{|\mathcal{O}|}{n}  $ to get
$t_{\IF}\left(\lambda_t/4; T_H(P),X_1^n\right)\le  4e^{-t}.$
}
\item There exists a constant $C_2 >0$ such that condition $\lambda_t \le \beta/2$ is verified for any $t\le C_2 n$.\\
\onepf{Use \ref{item:hyp_prop1:1} to get that $\lambda_t\le \beta/2$  is implied by
$$\frac{\beta}{4} \ge 6\sqrt{\frac{V_{\psi_H}}{n}}+8\sqrt{\frac{tv_{\psi_H}}{n}}+4\frac{Ct\beta}{n}.$$
Then, this is implied by the following system of equation with $C_1,C_2, C_3>0$ three numerical constants,
$$\begin{cases}
\beta \ge C_1\sqrt{\frac{V_{\psi_H}}{n}}\\
\beta \ge  \sqrt{C_2\frac{tv_{\psi_H}}{n}}\\
\beta \ge \frac{t\beta}{C_3 n}
\end{cases} \Leftarrow \begin{cases}
\beta \ge C_1\sqrt{\frac{V_{\psi_H}}{n}}\\
t\le C_2 n \beta^2/v_{\psi_H}\\
t \le C_3 n
\end{cases}
\Leftarrow \begin{cases}
\beta \ge C_1\sqrt{\frac{V_{\psi_H}}{n}}\\
t\le C_2 n \\
t \le C_3 n
\end{cases}$$
The last implication is from \ref{item:hyp_prop1:2}: $\beta^2 \ge 8 V_{\psi_H} \ge 8v_{\psi_H}$.
Inject this in Equation~\eqref{eq:tif_tt_h} to get the desired result.}
\end{npf}

\subsection{Proof of Proposition~\ref{prop:huber_corrupted_q}}
\begin{Hypopf}
\begin{enumerate}[label=(\roman{*}), ref=Hypothesis (\roman{*}), leftmargin=*]
\item $\E_P[\|X-\E_P[X]\|^q]<\infty$ for some $q\ge 2$. \label{hypopf:huber:qmoments}
\item $t \le C_1 n\min\left( \frac{\E_P[\|X-\E_P[X]\|^q]}{Tr(\Sigma)^{q/2}}, 1\right),$\label{hypopf:huber:t}
\end{enumerate}
\end{Hypopf}

\begin{npf}
\item $\beta = \left(\frac{2\E_P[\|X-\E_P[X]\|^q]}{8\varepsilon_n+C_3t/n}\right)^{1/q}$  verifies the conditions of Proposition~\ref{prop:huber_corrupted}.\\
\onepf{
We will check the stronger assumption that $\beta^2 \ge Tr(\Sigma)\max(8, C_1/n)$.\\
Having $\varepsilon_n \ge 0$, we can bound the chosen $\beta$ by
$\beta^q \ge n\frac{\E[\|X-\E[X]\|^q]}{C_3t}$, then there exists $C>0$ such that the hypothesis on $\beta$ is verified if
$$t \le C \frac{n\E[\|X-\E[X]\|^q]}{Tr(\Sigma)^{q/2}}.$$
This is condition on $t$ is verified in \ref{hypopf:huber:t}.}
\item  For any $t\le C_2 n $, we have \label{item:huber:1}
\begin{align*}
\|T_H(X_1^n)-\E_P[X]\|\le & 6\sqrt{\frac{V_{\psi_H}}{n}}+8\sqrt{\frac{tv_{\psi_H}}{n}}\\
&+4\left(\frac{C_3 t}{n}+ 8\varepsilon_n\right)^{1-1/q} \E[\|X-\E_P[X]\|^q]^{1/q}.
\end{align*}
\\
\onepf{
From the previous step, we can apply Proposition~\ref{prop:huber_corrupted}, and by Lemma~\ref{lem:bias_huber_using_concentration} and \ref{hypopf:huber:qmoments}, if $\beta^2 \ge V_{\psi_H}\max(8,C_1 /n)$ then with probability larger than $1-4e^{-t}-e^{-n/32}$,
\begin{align}\label{eq:huber2_1}
\|&T_H(X_1^n)-\E_P[X]\|\nonumber\\
&\le 6\sqrt{2\frac{V_{\psi_H}}{n}}+8\sqrt{2\frac{tv_{\psi_H}}{n}}+\frac{C_3t\beta}{n}+ 8\beta \varepsilon_n + \|\E_P[X]-T_H(P)\|\nonumber \\
&\le 6\sqrt{2\frac{V_{\psi_H}}{n}}+8\sqrt{2\frac{tv_{\psi_H}}{n}}+\frac{C_3t\beta}{n}+ 8\beta \varepsilon_n +\frac{2\E[\|X-\E[X]\|^q]}{(q-1)\beta^{q-1}}.
\end{align}
To simplify, we don't take into account the effect of $\beta$ on $V_{\psi_H}$ and $v_{\psi_H}$ when choosing $\beta$ then, we choose $\beta$ such that
$$\beta^q = \frac{2\E_P[\|X-\E_P[X]\|^q]}{8\varepsilon_n+C_3t/n}  $$
inject this in Equation~\eqref{eq:huber2_1} to get
\begin{multline*}
\|T_H(X_1^n)-\E_P[X]\|\le 6\sqrt{\frac{V_{\psi_H}}{n}}+8\sqrt{\frac{tv_{\psi_H}}{n}}\nonumber\\
+2^{1/q}\left( 1+\frac{1}{q-1}\right)\left(\frac{C_3 t}{n}+ 8\varepsilon_n\right)^{1-1/q} \E_P[\|X-\E_P[X]\|^q]^{1/q}.
\end{multline*}
The result follow because $2^{1/q}(1+\frac{1}{q-1}) \le 4$.
}
\item  For any $t\le C_2 n $,
\begin{multline*}\|T_H(X_1^n)-\E_P[X]\|\le 6\sqrt{2\frac{Tr(\Sigma)}{n}}\\
+8\sqrt{2\frac{t\|\Sigma\|_{op}}{n}}+C_4\left(\frac{C_3 t}{n}+8\varepsilon_n\right)^{1-1/q} \E[\|X-\E[X]\|^q]^{1/q}.
\end{multline*}\\
\onepf{
From Lemma~\ref{lem:bound_variance_variance}, $V_{\psi_H}\le Tr(\Sigma)$ and
\begin{align*}
v_{\psi_H} \le& \|\Sigma\|_{op}+ \|\E[X]-T_H(P)\|^2\\
\le& \|\Sigma\|_{op}+ \frac{4\E[\|X-\E[X]\|^{2q}]}{(q-1)^2\beta^{2q-2}}\le \|\Sigma\|_{op}\\
&+ \frac{2^{2/q}\E[\|X-\E[X]\|^q]^{2/q}}{(q-1)^2}\left(8\varepsilon_n+C_3\frac{t}{n}\right)^{2-2/q}
\end{align*}
by sub-linearity of the square root we obtain
$$\sqrt{v_{\psi_H}} \le  \sqrt{\|\Sigma\|_{op}}+\frac{2^{1/q}\E[\|X-\E[X]\|^q]^{1/q}}{(q-1)}\left(4\varepsilon_n+C_3\frac{t}{n}\right)^{1-1/q}$$
inject this in \ref{item:huber:1}, to get
\begin{align*}
\|&T_H(X_1^n)-\E_P[X]\|\\
\le& 6\sqrt{2\frac{Tr(\Sigma)}{n}}+8\sqrt{2\frac{t\|\Sigma\|_{op}}{n}}\\
&+\left(4 +2\sqrt{t/n}\right)\left(\frac{C_3 t}{n}+ 8\varepsilon_n\right)^{1-1/q} \E[\|X-\E[X]\|^q]^{1/q}\nonumber\\
\le& 6\sqrt{2\frac{Tr(\Sigma)}{n}}+8\sqrt{2\frac{t\|\Sigma\|_{op}}{n}}+C_4\left(\frac{C_3 t}{n}+8\varepsilon_n\right)^{1-1/q} \E[\|X-\E[X]\|^q]^{1/q}.
\end{align*}
}
\end{npf}

\subsection{Proof of Proposition~\ref{prop:huber_corrupted_sym}}
From Proposition~\ref{prop:huber_corrupted}, if $\beta^2 \ge V_{\psi_H}\max(8,C_1 /n)$ and having $\E_P[X]=T_H(P)$, we get with probability larger than $1-4e^-t-e^{-n/32}$,
\begin{align}\label{eq:hubersym_1}
\|T_H(X_1^n)-\E_P[X]\|&\le 6\sqrt{2\frac{V_{\psi_H}}{n}}+8\sqrt{2\frac{tv_{\psi_H}}{n}}+\frac{C_3t\beta}{n}+ 8\beta \varepsilon_n
\end{align}
Take $\beta = C_2\E[\|X-\E[X]\|]$ with $C_2^2 \ge \max(8,C_1)$, this value of $\beta$ verifies that $\beta^2 \ge V_{\psi_H}\max(8,C_1 /n)$ hence it verifies the condition of Proposition~\ref{prop:huber_corrupted}. Inject the value of $\beta$ in Equation~\eqref{eq:hubersym_1} to get,
\begin{align}
\|T_H(X_1^n)-\E_P[X]\|&\le 6\sqrt{\frac{V_{\psi_H}}{n}}+8\sqrt{\frac{tv_{\psi_H}}{n}}+C_2\E[\|X-T(P)\|]\left(\frac{C_3t}{n}+ 8 \varepsilon_n\right).
\end{align}
\subsection{Proof of Proposition~\ref{prop:catoni_corrupted_q}}\label{sec:proof_catoni_corrupted_2}
\begin{Hypopf}
\begin{enumerate}[label=(\roman{*}), ref=Hypothesis (\roman{*}), leftmargin=*]
\item With probability larger than $1-\delta_O$, we have \\
$ \frac{1}{|\mathcal{O}|}\sum_{i\in \mathcal{O}} \psi_C(\|X_i-T_C(P)\|)+\psi_C(\|X'_i-T_C(P)\|) \le \beta C_O.$
\item $\beta^3=\frac{5\E[\|X-\E[X]\|^3]}{120C_O\varepsilon_n+24Ct\log(n)/n},$
\item $|\mathcal{O}|\le n/(20C_O)$
\end{enumerate}
\end{Hypopf}
\begin{npf}
\item $\beta$ verify that  $\beta^2 \ge 10Tr(\Sigma)$. \\
\onepf{Having $\varepsilon_n \ge 0$ and $s\ge 1$, we can bound the chosen $\beta$ by
$\beta^3 \ge \frac{5\E[\|X-\E[X]\|^3]}{48 C t\log(n)/n}$. Then, because there exists $C_3>0$ such that $t \le C_3n / \log(n)$, we have $\beta^2 \ge 10\E[\|X-\E[X]\|^3]^{2/3}$ (if needed, $C_3$ can be decreased to obtain  the right constant), hence $\beta^2 \ge 10Tr(\Sigma)$. }

\item \begin{multline*}
t_{\IF}(\lambda; T_C(P), X_1^n)\\\le  \P\left(\left\|\frac{1}{n}\sum_{i =1}^n \frac{X_i'-T_C(P)}{\|X_i'-T_C(P)\|}\psi_C\left(\|X_i'-T_C(P)\|\right) \right\| \ge \lambda /5 - \varepsilon_n C_O\beta \right)+ \delta_O
\end{multline*}
 \label{item:catoni:3}\\
\onepf{
From Theorem~\ref{th:tIF_tT_1} with $\gamma=4/5$, we have that if $V_{\psi_C}=\E_P[\psi_C(\|X-T_C(P)\|)^2]\le \psi_C(\beta/2)^2/2=\beta^2 \log(13/8)^2/2$, then for all $\lambda \in (0,\beta/2)$,
\begin{equation}\label{eq:tif_tt_c}
t_T(\lambda; T_C(P),X_1^n)\le t_{\IF}\left(\lambda/5; T_C(P),X_1^n\right)+e^{-n/50}.
\end{equation}
The condition $V_{\psi_C}\le \beta^2 \log(13/8)^2/2$ is verified because $\beta^2 \ge 10 Tr(\Sigma)$ and $V_{\psi_C}\le Tr(\Sigma)$ by Lemma~\ref{lem:bound_variance_variance}.

Now, we take care of the outliers in the right-hand side of Equation~\eqref{eq:tif_tt_c},
\begin{align}\label{eq:tif_c_cor}
t_{\IF}&\left(\lambda/5; T_C(P)\right)\nonumber\\
\le& \P\left(\left\|\frac{1}{n}\sum_{i=1}^n \frac{X_i'-T_C(P)}{\|X_i'-T_C(P)\|}\psi_C\left(\|X_i'-T_C(P)\|\right) \right\| \ge \lambda /5 \right. \nonumber\\
&\left. - \frac{1}{n}\sum_{i\in \mathcal{O}} \psi_C(\|X_i-T_C(P)\|)+\psi_C(\|X'_i-T_C(P)\|) \right)\nonumber\\
\le&  \P\left(\left\|\frac{1}{n}\sum_{i =1}^n \frac{X_i'-T_C(P)}{\|X_i'-T_C(P)\|}\psi_C\left(\|X_i'-T_C(P)\|\right) \right\| \ge \lambda /5 - \varepsilon_n C_O\beta \right)+ \delta_O
\end{align}
}
\item
\label{eq:catoni:bias}
$\E[\|X- T_C(P)\|^2] \le 2Tr(\Sigma).$\\
\onepf{
We apply Lemma~\ref{lem:bias_general} with $\gamma=4/5$, $\|\psi'\|_{\infty}\le 1$ and two finite moments to get, $\|T_C(P)-\E_P[X]\| \le \frac{5Tr(\Sigma)}{4\beta}$. We get $\|T_C(P)-\E_P[X]\| \le 5Tr(\Sigma)/(4\sqrt{10}) \le Tr(\Sigma)$.
Hence, by bias-variance decomposition, $\E[\|X- T_C(P)\|^2] = \E[\|X-\E[X]\|^2] +\|T_C(P)-\E_P[X]\|^2 \le 2Tr(\Sigma).$

}
\item We have \label{item:catoni:4}
$$t_{\IF}\left(\lambda_t/5; T_C(P)\right)\le 4e^{-t}+\delta_0$$
where $$\lambda_t = 5\left( \frac{3}{2}\sqrt{\frac{V_{\psi_C}}{n}}+2\sqrt{\frac{tv_{\psi_C}}{n}}+\frac{Ct\beta}{n}\log(n) \right)+ 5\varepsilon_n C_O\beta.$$\\
\onepf{
In order to use Lemma~\ref{lem:subgaussian}, we use the following Lemma.
\begin{Lemma}\label{lem:catoni_moments}
If $X$ satisfies $\E_P[\|X\|^2]<\infty$, then, for all $q\in \N^*$,
\[
\E_P[\psi_C(\|X-T_C(P)\|)^{q}]\le  q!(s\beta)^{q},
\]
where
$$s=\max\left(e,\log\left(1+\frac{\E_P[\|X-T_C(P)\|])}{\beta}+\frac{\E[\|X-T_C(P)\|^2]}{2\beta^2} \right)\right).$$
\end{Lemma}

The proof of Lemma~\ref{lem:catoni_moments} is postponed to Section~\ref{sec:proof_catoni_moments}. By \ref{eq:catoni:bias} and the fact that $\beta^2 \ge Tr(\Sigma)$, we have $s \le \max(e, \log(1+1+1/2))\le e$. Then, by Lemma~\ref{lem:catoni_moments}, for any $q \in \N^*$,
\[
\E_P[\psi_C(\|X-T_C(P)\|)^{q}]\le  q!(e\beta)^{q},
\]

Then, using the power series expansion of the exponential function, we get that, for all $t>\beta e$,
\begin{align*}
\E\left[\exp\left(\frac{\psi_C(\|X-T_C(P)\|)}{t} \right)\right]&=\sum_{q=0}^\infty \frac{\E[\psi_C(\|X-T(P)\|^q]}{t^q q!}\\
&\le \sum_{q=0}^{\infty}\frac{\beta^qe^q}{t^q}=\frac{1}{1-\beta e/t}.
\end{align*}
Choosing $t=2\beta e$ shows that $\|\psi_C(\|X-T_C(P)\|)\|_{\psi_1}\le 2e\beta $. Then, using Lemma~\ref{lem:subgaussian}, we get for all $t>0$, with probability larger than $1-4e^{-t}$,
$$\left\|\frac{1}{n}\sum_{i=1}^n \frac{X_i'-T_C(P)}{\|X_i'-T_C(P)\|}\psi_H\left(\|X_i'-T_C(P)\|\right) \right\| \le  \frac{3}{2}\sqrt{\frac{V_{\psi_C}}{n}}+2\sqrt{\frac{tv_{\psi_C}}{n}}+\frac{Ct\beta}{n}\log(n) ,$$
Hence, from \ref{item:catoni:3}, we have
$$t_{\IF}\left(\lambda_t/5; T_C(P)\right)\le 4e^{-t}+\delta_0$$
where $$\lambda_t = 5\left( \frac{3}{2}\sqrt{\frac{V_{\psi_C}}{n}}+2\sqrt{\frac{tv_{\psi_C}}{n}}+\frac{Ct\beta}{n}\log(n) \right)+ 5\varepsilon_n C_O\beta.$$}
\item The condition $\lambda_t \le \beta/2$ is implied by $t \le C_3 n /\log(n)$ for some absolute constant $C_3$.\\
\onepf{
We use that $\varepsilon_nC_O \le 1/20$ to get that $\lambda_t\le \beta/2$ is implied by
$$\frac{\beta}{4} \ge \frac{15}{2}\sqrt{\frac{V_{\psi_C}}{n}}+10\sqrt{\frac{tv_{\psi_C}}{n}}+5\frac{Cts\beta}{n}\log(n).$$
Then, this is implied by the following system of equation with $C_1,C_2, C_3>0$ three numerical constants,
$$\begin{cases}
\beta \ge C_1\sqrt{\frac{V_{\psi_C}}{n}}\\
\beta \ge  \sqrt{C_2\frac{tv_{\psi_C}}{n}}\\
\beta \ge \frac{t \beta}{C_3 n}\log(n)
\end{cases} \Leftarrow \begin{cases}
\beta \ge C_1\sqrt{\frac{V_{\psi_C}}{n}}\\
t\le C_2 n \beta^2/v_{\psi_C}\\
t \le C_3 \frac{n}{s\log(n)}
\end{cases}
\Leftarrow \begin{cases}
\beta \ge C_1\sqrt{\frac{V_{\psi_C}}{n}}\\
t\le C_2 n \\
t \le C_3 \frac{n}{\log(n)}
\end{cases}$$
The last implication is because $\beta^2 \ge 8 V_{\psi_C} \ge 8v_{\psi_C}$. The first inequality is necessarily verified because $\beta^2 \ge 10Tr(\Sigma)$, hence the only remaining condition is $t \le C_3 n /\log(n)$ for some absolute constant $C_3$.
}
\item There exists a numerical constant $C_2>0$ such that with probability larger than $1-4e^{-t}-\delta_0$,
\begin{multline}\label{eq:catoni_corrupt_2}
\|T_C(X_1^n)-\E_P[X]\|\le \frac{15}{2}\sqrt{\frac{Tr(\Sigma)}{n}}+10\sqrt{\frac{t\|\Sigma\|_{op}}{n}}\\+C_2\E[\|X-\E[X]\|^3]^{1/3}\left(C_O\varepsilon_n+Ct\frac{\log(n)}{n}\right)^{2/3}.
\end{multline}
\onepf{
From Lemma~\ref{lem:bias_general}, with $3$ finite moments, we have
\begin{equation}\label{eq:bound_bias_catoni}
\|T_C(P)-\E_P[X]\| \le \frac{5\E[\|X-\E[X]\|^3]}{24\beta^2}.
\end{equation}
Then, from \ref{item:catoni:4} with probability larger than $1-4e^{-t}-e^{-n/50}$,
\begin{multline*}\|T_C(X_1^n)-\E_P[X]\|\le \frac{15}{2}\sqrt{\frac{V_{\psi_C}}{n}}\\
+10\sqrt{\frac{tv_{\psi_C}}{n}}+\frac{Ct\beta}{n}\log(n)+ 5 \beta C_O\varepsilon_n+ \frac{5\E[\|X-\E[X]\|^3]}{24\beta^2}.
\end{multline*}

Then, there exists a constant $C_4>0$ such that
\begin{multline}\label{eq:catoni2_2}
\|T_C(X_1^n)-\E_P[X]\|\le \frac{15}{2}\sqrt{\frac{V_{\psi_C}}{n}}+10\sqrt{\frac{tv_{\psi_C}}{n}}\\
+C_4\E[\|X-\E[X]\|^3]^{1/3}\left(C_O\varepsilon_n+t\frac{\log(n)}{n}\right)^{2/3},
\end{multline}
To simplify, we also use the bounds on the variance from~\ref{lem:bound_variance_variance}, $V_{\psi_C}\le Tr(\Sigma)$ and also
\begin{align*}
v_{\psi_C} &\le \|\Sigma\|_{op}+ \|\E[X]-T_C(P)\|^2\\
&\le \|\Sigma\|_{op}+ \frac{25\E[\|X-\E[X]\|^3]^2}{576\beta^2}.
\end{align*}
By sub-linearity of the square root we obtain
\begin{align*}\sqrt{v_{\psi_C}} &\le  \sqrt{\|\Sigma\|_{op}}+\frac{5\E[\|X-\E[X]\|^3]}{24\beta^2}\\
&= \sqrt{\|\Sigma\|_{op}}+C_4\E[\|X-\E[X]\|^3]^{1/3}\left(C_O\varepsilon_n+ts\frac{\log(n)}{n}\right)^{2/3},
\end{align*}
inject this in Equation~\eqref{eq:catoni2_2}, to get
\begin{multline*}
\|T_C(X_1^n)-\E_P[X]\|\le \frac{15}{2}\sqrt{\frac{Tr(\Sigma)}{n}}+10\sqrt{\frac{t|\Sigma\|_{op}}{n}}\\
+C_4\left(1+ \sqrt{\frac{2t}{n}}\right)\E[\|X-\E[X]\|^3]^{1/3}\left(C_O\varepsilon_n+Ct\frac{\log(n)}{n}\right)^{2/3}.
\end{multline*}
Having $t\lesssim n$, we get that there exists a numerical constant $C_2>0$ such that
\begin{multline}\label{eq:catoni_corrupt_2}
\|T_C(X_1^n)-\E_P[X]\|\le \frac{15}{2}\sqrt{\frac{Tr(\Sigma)}{n}}+10\sqrt{\frac{t\|\Sigma\|_{op}}{n}}\\+C_2\E[\|X-\E[X]\|^3]^{1/3}\left(C_O\varepsilon_n+Ct\frac{\log(n)}{n}\right)^{2/3}.
\end{multline}}
\item Using the sub-linearity of $x \mapsto x^{2/3}$ and using that 
$$t \le C_2  \frac{n}{\log(n)^4}\frac{\sqrt{\|\Sigma\|_{op}}}{\E[\|X-\E[X]\|^3]^{1/3}}$$ 
for $C_2$ small enough in equation~\eqref{eq:catoni_corrupt_2}, we have 
\begin{multline*}C_2\E[\|X-\E[X]\|^3]^{1/3}\left(C_O\varepsilon_n+Ct\frac{\log(n)}{n}\right)^{2/3}\\
\le C_3\E[\|X-\E[X]\|^3]^{1/3}\varepsilon_n^{2/3}+\sqrt{\frac{t\|\Sigma\|_{op}}{n}}
\end{multline*}
which finishes the proof.
\end{npf}

%
 \subsection{Proof of Proposition~\ref{prop:poly_corrupted_2}}\label{sec:proof_poly2}
\begin{Hypopf}
\begin{enumerate}[label=(\roman{*}), ref=Hypothesis (\roman{*}), leftmargin=*]
\item With probability larger than $1-\delta_O$, we have \\
$ \frac{1}{|\mathcal{O}|}\sum_{i\in \mathcal{O}} \psi_P(\|X_i-T_P(P)\|)+\psi_P(\|X'_i-T_P(P)\|) \le \beta C_O.$\label{hypopf:poly:1}
\item $\beta^2 = Tr(\Sigma)$\label{hypopf:poly:2}
\end{enumerate}
\end{Hypopf}
\begin{npf}
\item For $t \le C_1 n$, with probability larger than $1-e^{-C_1t}-e^{-n/512}-\delta_O$,
\begin{align*}
\|T_P(X_1^n)-T_P(P)\|\le 16 \sqrt{\frac{\E_P[\|X-T_P(P)\|^2]t }{n}}+ 16\varepsilon_n C_O\beta.
\end{align*}
\onepf{From Theorem~\ref{th:tIF_tT_1} with $\gamma=1/4$, we have that if $V_{\psi_P}=\E_P[\psi_P(\|X-T_P(P)\|)^2]\le \psi_P(\beta/2)^2/2=\beta^2/8$, then for all $\lambda \in (0,\beta/2)$,
\begin{equation}\label{eq:tif_tt_p}
t_T(\lambda; T_P(P),X_1^n)\le t_{\IF}\left(\lambda/16; T_P(P),X_1^n\right)+e^{-n/512}.
\end{equation}

The following lemma applies.
\begin{Lemma}\label{lem:concentration_poly}
Let $n \in \N^*$, suppose $X_1,\dots,X_n$ are i.i.d.
Let $q \in \N^*$ and suppose $\E_P[\|X\|^{q}]<\infty$.
There exists an absolute constant $K>0$ such that
$$t_{\IF}(\lambda)\le \frac{\E_P[\|X-T_P(P)\|^q]}{\beta^q}\left(\frac{Kpq\beta}{\sqrt{n}\lambda}\right)^{qp}.$$
\end{Lemma}
The proof is postponed to Section~\ref{sec:proof_concentration_poly}.\\
 Take $q=2$ and $\lambda = \sqrt{\E_P[\|X-T_P(P)\|^2]t/n}$, we get for all $t>0$,
\begin{align*}\P&\left(\left\|\frac{1}{n}\sum_{i=1}^n\frac{X_i'-T_P(P)}{\|X_i'-T_P(P)\|}\psi_P(\|X_i'-T_P(P)\|) \right\| \ge\sqrt{\frac{\E_P[\|X-T_P(P)\|^2]t }{n}}\right)\\
&\le \left(\frac{\beta^2}{\E_P[\|X-T_P(P)\|^2]}\right)^{p-1}\left(\frac{4K^2p}{t}\right)^{2p}.
\end{align*}
Take having $Tr(\Sigma)\le \E_P[\|X-T(P)\|^2]$, we have $\beta^2 \le \E_P[\|X-T(P)\|^2]$ hence
\begin{align*}
\P&\left(\left\|\frac{1}{n}\sum_{i=1}^n\frac{X_i'-T_P(P)}{\|X_i'-T_P(P)\|}\psi_P(\|X_i'-T_P(P)\|) \right\| \ge\sqrt{\frac{\E_P[\|X-T_P(P)\|^2]t }{n}}\right)\\
&\le\left(\frac{4K^2p}{t}\right)^{2p}.
\end{align*}
Take $p=\frac{t}{4K^2}e^{-1}$, and there is a constant $C_1>0$ such that with probability larger than $1-e^{-C_1t}$, we have
\begin{equation}\label{eq:tif_poly}
\left\|\frac{1}{n}\sum_{i=1}^n\frac{X_i'-T_P(P)}{\|X_i'-T_P(P)\|}\psi_P(\|X_i'-T_P(P)\|) \right\| \le\sqrt{\frac{\E_P[\|X-T_P(P)\|^2]t }{n}}.
\end{equation}

Now, we take care of the outliers in a similar way as in Proposition \ref{prop:huber_corrupted} and Proposition~\ref{prop:catoni_corrupted_q}. We have from \ref{hypopf:poly:1},
\begin{align}\label{eq:tif_c_cor}
t_{\IF}&\left(\lambda/16; T_P(P),X_1^n\right)\le  t_{\IF}\left( \frac{\lambda}{16}- \varepsilon_n C_O\beta; T_P(P), (X')_1^n \right)+ \delta_O
\end{align}
Hence, from Equation~\eqref{eq:tif_poly} and Equation~\eqref{eq:tif_c_cor} with 
$$\lambda=\lambda_t=16\sqrt{\frac{\E_P[\|X-T_P(P)\|^2]t }{n}}+ 16\varepsilon_n C_O\beta,$$ we have with probability larger than $1-e^{-C_1t}-e^{-n/512}-\delta_O$,
\begin{align*}
\|T_P(X_1^n)-T_P(P)\|&\le 16 \sqrt{\frac{\E_P[\|X-T_P(P)\|^2]t }{n}}+ 16\varepsilon_n C_O\beta \\
&\le 16 \sqrt{\frac{\E_P[\|X-T_P(P)\|^2]t }{n}}+ 16\varepsilon_n C_O\beta.
\end{align*}}
\item There exist a constant $C_1'>0$ such that for any $t \le C_1' n $, we have $\lambda_t \le \beta/2$.\\
\onepf{The condition $\lambda_t \le \beta/2$ is implied by
$$
 16 \sqrt{\frac{\E_P[\|X-T_P(P)\|^2]t }{n}}\le \beta/4$$
 using the fact that $\varepsilon_n \le 1/(64C_O)$. This simplifies with $t \lesssim n\beta^2/\E_P[\|X-T_P(P)\|^2]$ and then, by \ref{hypopf:poly:2} $t \lesssim n$. }
 \item With probability larger than $1-e^{-C_2t}-e^{-n/512}-\delta_O$, we have
$$\|T_P(X_1^n)-\E[X]\|\le 16 \sqrt{65\frac{Tr(\Sigma)t }{n}}+ 16\varepsilon_n C_O\sqrt{ Tr(\Sigma)}.
$$\\
\onepf{
For any $t \le C_1 n\beta^2/\E[\|X-T_P(P)\|^2]$, with probability larger than $1-e^{-C_2t}-e^{-n/512}-\delta_O$,
\begin{equation}\label{eq:poly_1}
\|T_P(X_1^n)-T_P(P)\|\le 16 \sqrt{\frac{\E[\|X-T_P(P)\|^2]t }{n}}+ 16\varepsilon_n C_O\beta.
\end{equation}
From Lemma~\ref{lem:bias_general}, we also have that $\|T_P(P)-\E[X]\| \le 8Tr(\Sigma)/\beta $. Hence,
$$\E[\|X-T_P(P)\|^2]\le \E[\|X-\E[X]\|^2]+\|\E[X]-T_P(P)\|^2\le Tr(\Sigma)\left(1+\frac{64Tr(\Sigma)}{\beta^2}\right)$$
from \ref{hypopf:poly:2} we have $\beta^2 = Tr(\Sigma)$, hence
$$\E[\|X-T_P(P)\|^2]\le Tr(\Sigma)(1+64)=65Tr(\Sigma),$$
and $\|T_P(P)-\E[X]\|\le 8\sqrt{Tr(\Sigma)}$. Then, inject this in Equation~\eqref{eq:poly_1},
$$\|T_P(X_1^n)-\E[X]\|\le 16 \sqrt{65\frac{Tr(\Sigma)t }{n}}+  16\varepsilon_n C_O\sqrt{ Tr(\Sigma)}.
$$
}
\end{npf}
\section{Proof of Lemmas}
\subsection{Proof of Lemma~\ref{lem:convexity}}\label{sec:proof_convexity}
By direct differentiation of  $Z_{P,\psi}$,
\begin{align*}
Jac(Z_{P,\psi})(\theta)=&\E\left[\frac{(X-\theta)\otimes(X-\theta)}{\|X-\theta\|^3} \psi\left(\left\|X-\theta\right\|\right) \right]-\E\left[\frac{I_\mathcal{H}}{\|X-\theta\|} \psi\left(\left\|X-\theta\right\|\right) \right]\\
&-\E\left[\frac{(X-\theta)\otimes(X-\theta)}{\|X-\theta\|^2} \psi'\left(\left\|X-\theta\right\|\right) \right]
\end{align*}
Hence, for any $u \in S$,
\begin{multline}\label{eq:diff_Z}
u^TJac(Z_{P,\psi})(\theta)u = -\E\left[\frac{1}{\|X-\theta\|} \psi\left(\left\|X-\theta\right\|\right) \right]+\E\left[\frac{\langle X-\theta, u \rangle^2}{\|X-\theta\|^3} \psi\left(\left\|X-\theta\right\|\right) \right]\\
-\E\left[\frac{\langle X-\theta, u \rangle^2}{\|X-\theta\|^2} \psi'\left(\left\|X-\theta\right\|\right) \right]
\end{multline}
Then, remark that because $\psi$ is concave, by Taylor's theorem we have $\forall y \ge 0$, $\psi(y)\ge \psi(0)+y\psi'(y)$. And then, because and $\psi(0)=0$ we have for any $y\ge 0,$ $\psi(y)\ge y\psi'(y)$. \\
Moreover, by Cauchy-Schwarz inequality we also have
$1-\left\langle \frac{X-\theta}{\|X-\theta\|},u\right\rangle^2\ge 0$.
Hence, with these two inequalities, we get
$$\E\left[\frac{1-\langle \frac{X-\theta}{\|X-\theta\|},u\rangle^2 }{\|X-\theta\|} \psi\left(\left\|X-\theta\right\|\right) \right]\ge \E\left[\left(1-\langle \frac{X-\theta}{\|X-\theta\|},u\rangle^2 \right) \psi'\left(\left\|X-\theta\right\|\right) \right].$$
Inject this in \eqref{eq:diff_Z} to get the result: for any  $u  \in S$,
 $$u^TJac(Z_{P,\psi})(\theta)u\le -\E\left[\psi'\left(\left\|X-\theta\right\|\right) \right].$$
\subsection{Proof of Lemma~\ref{lem:unicity}}\label{sec:proof_unicity}

\begin{Hypopf}
\begin{enumerate}[label=(\roman{*}), ref=(\roman{*}), leftmargin=*]
\item $\E[\rho(\|X-\E[X]\|)]<\rho(\beta)$ \label{hypopf:1}
\end{enumerate}
\end{Hypopf}
\begin{npf}
\item $T(P)$ is a critical point of $J(\theta)=\E[\rho(\|X-\theta\|)]$.\\
\onepf{First, notice that from Assumptions \ref{ass:1} and \ref{ass:3}
$\rho$ is two times derivable and increasing on $\R_+$. Hence, $J$ is differentiable and its gradient is
$$\nabla J(\theta)=-\E\left[\frac{X-\theta}{\|X-\theta\|}\psi(\|X-\theta\|)\right].$$
Then, by definition of $T(P)$ (Equation~\eqref{eq:def_T_th}), $T(P)$ verifies $\nabla J(T(P))=0$, i.e. $T(P)$ critical point of $J$.
}

\item \label{claim:convex}{$J$ is convex }\\
\onepf{
 Let $\mathrm{Hess}(J)$ denote the Hessian of $J$.
From Lemma~\ref{lem:convexity}, for any $\theta \in \mathcal{H}$  and $u \in S$,\label{item:hessian3}
\begin{align*}
u^T\mathrm{Hess}(J)(\theta)u& = -u^TJac(Z_\psi)(\theta)u\ge  \E\left[\psi'(\|X-\theta\|)\right]\ge 0
\end{align*}
Hence, $J$ is convex because its Hessian is positive.}
\item  $T(P)$ exists\label{item:minimum}\\
\onepf{Because $\rho$ is increasing, we have $\rho(x)\xrightarrow[x \to\infty]{}\infty$ and hence $J(\theta)\xrightarrow[\|\theta\| \to \infty]{}\infty$. Hence, $J$ is coercive and as $J$ is also convex (\ref{claim:convex}), its minimum $T(P)$ exists.
}
\item \label{claim:strict_convex}{$J$ is strictly convex at $T(P)$.}\\
\onepf{
From  Assumption~\ref{ass:3}, and because $\rho$ is increasing,
\begin{align*}
\E[\psi'(\|X-T(P)\|)]&\ge \gamma \P(\|X-T(P)\|\le \beta)= \gamma \P(\rho(\|X-T(P)\|)\le \rho(\beta)).
\end{align*}
Hence, by Markov's inequality,
\begin{equation}\label{eq:item4}
\E[\psi'(\|X-T(P)\|)]\ge  \gamma \left(1-\frac{\E[\rho(\|X-T(P)\|)]}{\rho(\beta)}\right).
\end{equation}
By \ref{item:minimum}, we have that $T(P)$ minimizer of $J$, hence
$$\frac{\E[\rho(\|X-T(P)\|)]}{\rho(\beta)}= \frac{J(T(P))}{\rho(\beta)}\le \frac{J(\E[X])}{\rho(\beta)}.$$
Inject this in Equation~\eqref{eq:item4} to get
$$\E[\psi'(\|X-T(P)\|)]\ge  \gamma \left( 1-\frac{\E[\rho(\|X-\E[X])\|)]}{\rho(\beta)}\right).$$
Then, using Hypothesis \ref{hypopf:1}, we get $\E[\psi'(\|X-T(P)\|)]>0$, which implies from Lemma~\ref{lem:convexity} that for all $u\in\mathcal{H}$, $u\neq 0$,
$u^T\mathrm{Hess}(J)(T(P))u>0$. Hence $J$ is strictly convex at $T(P)$, the minimizer of $J$.
}
\item $T(P)$ as defined by Equation~\eqref{eq:def_T_th} exists and is unique.\\
\onepf{As $T(P)$ is minimizer of $J$, it is also root of Equation~\eqref{eq:def_T_th} and the existence and unicity we proven in \ref{claim:convex} and \ref{claim:strict_convex}. }
\end{npf}

\subsection{Proof of Lemma~\ref{lem:ass_exemples}}\label{sec:proof_ass_exemples}
\noindent\textbf{Huber's score function:} The equality for the Huber's score function is immediate by derivation of $\psi_H$.

\noindent\textbf{Catoni's score function:} $\psi_C$ is differentiable, and we have for all $x\ge 0$,
$$\psi_{C}'(x)=\frac{1+\frac{x}{\beta}}{1+ \frac{x}{\beta}+\frac{x^2}{2\beta^2}}.$$
This function is decreasing on $\R_+$, positive and even, hence
$$\psi_{C}'(x)\ge \psi_{C}'(\beta)\1\{x\le \beta\}= \frac{4}{5}\1\{x\le \beta \}.$$

\noindent\textbf{Polynomial score function:} $\psi_P$ is differentiable, and we have for all $x\ge 0$
$$\psi_P(x)'=\frac{1+\frac{1}{p}\left(\frac{x}{\beta}\right)^{1-1/p}}{\left(1+\left(\frac{x}{\beta}\right)^{1-1/p}\right)^2} .$$
As in the case of Catoni's score function, this function is decreasing over $\R_+$, positive and even. Then, we get
$$\psi_P'(x)\ge \psi_P'(\beta)\1\{x\le \beta\}= \frac{1}{4}\left(1+\frac{1}{p}\right)\1\{x\le \beta\}\ge \frac{1}{4}\1\{x\le \beta\}. $$

\subsection{Proof of Lemma~\ref{lem:bias_huber_using_concentration}}\label{sec:proof_bias_huber}
From Theorem~\ref{th:cornerstone} we only need to control $Z_\beta(\E[X])$. We have,
\begin{align*}
Z_\beta&(\E[X])\\
&=\E\left[\beta\frac{(X-\E[X])}{\|X-\E[X]\|} \psi_1\left(\left\|\frac{X-\E[X]}{\beta}\right\|\right) \right]\\
&=\E\left[\beta\frac{(X-\E[X])}{\|X-\E[X]\|} \psi_1\left(\left\|\frac{X-\E[X]}{\beta}\right\|\right) \right]-\E\left[\beta\frac{(X-\E[X])}{\|X-\E[X]\|} \left\|\frac{X-\E[X]}{\beta}\right\| \right]
\end{align*}
Hence, by triangular inequality,
\begin{align*}
\|Z_\beta(\E[X])\|&\le \beta\E\left[\left|\psi_1\left(\left\|\frac{X-\E[X]}{\beta}\right\|\right) - \left\|\frac{X-\E[X]}{\beta}\right\|  \right|\right]
\end{align*}
We denote $Y=\|X-\E[X]\|/\beta$, we have
\begin{align*}
\|Z_\beta(\E[X])\|\le \beta\E\left[\left|\psi_1\left(Y\right) -Y\right|\right]&=\beta\int |\psi_1(y)-y|dF_Y(y)\\
&=\beta\int_{0}^{\infty} (y-\psi_1(y))dF_Y(y)
\end{align*}
because $\psi_1$ is $1$-Lipshitz and $\psi_1(0)=0$. Then, by integration by part,
$$\|Z_\beta(\E[X])\|\le \beta\int_{0}^{\infty} (1-\psi_1'(y))(1-F_Y(y))dy
$$
Until now, the proof was valid for any $\psi_1$, for the specific case of Huber score function and using Theorem~\ref{th:cornerstone}, we get that
$$ \|\E[X]-T_H(P)\|\le  2\beta\int_{1}^{\infty}\P\left(\|X-\E[X]\| \ge \beta y \right)dy$$
Then, use Markov's inequality,
$$ \|\E[X]-T_H(P)\|\le  2\beta\int_{1}^{\infty}\frac{\|X-\E[X]\|^q}{\beta^q y^q}dy=\frac{2\|X-\E[X]\|^q}{(q-1)\beta^{q-1}}.$$

\subsection{Proof of Lemma~\ref{lem:bias_general}}\label{sec:proof_bias_general}
We have,
\begin{align*}
Z_\beta(\E[X])&=\E\left[\frac{X-\E[X]}{\|X-\E[X]\|}\beta\psi_1\left(\left\|\frac{X-\E[X]}{\beta}\right\|\right)\right]
\end{align*}
Then, by Taylor expansion
\begin{align*}
\left\|Z_\beta(\E[X])\right\|\le & \left\|\E\left[\frac{X-\E[X]}{\|X-\E[X]\|}\beta\left\|\frac{X-\E[X]}{\beta}\right\|\right]\right\|+\beta\E\left[\frac{\|\psi_1^{(k)}\|_\infty}{k!}\left\|\frac{X-\E[X]}{\beta}\right\|^k\right]\\
=&\E\left[\frac{\|\psi_1^{(k)}\|_\infty\E[\|X-\E[X]\|^k]}{k!\beta^{k-1}}\right]
\end{align*}
which proves the first part of the lemma using Theorem~\ref{th:cornerstone}.
In the case of Bernoulli distribution, the result follows from a Taylor expansion:
\begin{align*}
Z_\beta(\E[X])=&\E\left[\sign(X-\E[X])\beta\psi_1\left(\left|\frac{X-\E[X]}{\beta}\right|\right)\right]\\
=&p\left(\beta\psi_1\left(\frac{1-p}{\beta}\right)\right)+(1-p)\left(\beta\psi_1\left(\frac{-p}{\beta}\right)\right)\\
=&p\beta\left(\frac{1-p}{\beta}+\frac{1}{k!}\psi_1^{(k)}(0)\frac{(1-p)^k}{\beta^k}+o\left(\frac{1}{\beta^k}\right)\right)\\
&-(1-p)\beta\left(\frac{p}{\beta}+\frac{1}{k!}\psi_1^{(k)}(0)\frac{p^k}{\beta^k}+o\left(\frac{1}{\beta^k} \right)\right).
\end{align*}
\subsection{Proof of Lemma~\ref{lem:bound_variance_variance}}\label{sec:proof_bound_variance}
\begin{npf}
\item For any $\psi$ that satisfies Assumptions~\ref{ass:prop_psi}, $V_\psi\le \mathrm{Tr}(\Sigma)$.\\
\onepf{First, remark that we have for all $x \in \R_+$, $\psi_1^2(x)\le 2\rho_1(x)$. Indeed, let $h(x)=\psi_1^2(x)-2\rho_1(x)$, its derivative is $h'(x)=2\psi_1(x)\left(\psi_1'(x)-1\right)$ and because  $\psi_1'\le 1$ and $\psi_1(0)=0$, we get that $h$ is decreasing, the fact that $h(0)=0$ implies that for all $x \in \R_+$, $\psi_1^2(x)\le 2\rho_1(x)$. Then,
\begin{equation}\label{eq:bound_V1}
V_\psi\le 2\beta^2\E\left[\rho_1\left(\frac{\|X-T(P)\|}{\beta}\right)\right].
\end{equation}
Define $J(\theta)=\E\left[\rho_1\left(\frac{\|X-\theta\|}{\beta}\right)\right]$ by definition~\ref{def:T_min}, $T(P)$ is the minimum of $J$ and by equation~\eqref{eq:bound_V1},
$$V_\psi\le 2\beta^2\E\left[\rho_1\left(\frac{\|X-\E[X]\|}{\beta}\right)\right]$$
Then finally, using that by integration of $\psi_1'\le 1$ we have $\rho_1(x)\le x^2/2$, hence the result. }

\item For any $\psi$ that satisfies Assumptions~\ref{ass:prop_psi}, $x_\psi\le \|\Sigma\|_{op}+\|\E[X]-T(P)\|^2$.\\
\onepf{Note that because $\psi_1(x)\le 2\rho_1(x)\le x^2$,
\begin{align*}
v_\psi &= \beta^2\sup_{u \in S} \E\left[\frac{\langle u, X-T(P)\rangle ^2}{\|X-T(P)\|^2}\psi_1\left(\frac{\|X-T(P)\|}{\beta}\right)^2  \right]\\
&\le \beta^2\sup_{u \in S} \E\left[\frac{\langle u, X-T(P)\rangle ^2}{\|X-T(P)\|^2}\left(\frac{\|X-T(P)\|}{\beta}\right)^2 \right]\\
&=\sup_{u \in S} \E\left[\langle u, X-T(P)\rangle ^2 \right]\\
&=\sup_{u \in S} \E\left[\left(\langle u, X-\E[X]\rangle + \langle u,\E[X]-T(P)\rangle\right) ^2 \right]\\
&=\sup_{u \in S} \E\left[\langle u, X-\E[X]\rangle^2 + \langle u,\E[X]-T(P)\rangle^2\right]=\|\Sigma\|_{op}+\|\E[X]-T(P)\|^2
\end{align*}}
\end{npf}
\subsection{Proof of Lemma~\ref{lem:lower_bound_variance_huber}}\label{sec:proof_lower_variance}

\begin{npf}
\item if $X$ has $q>2$ finite moments, then \\
\begin{multline}\label{eq:var_psih}
V_{\psi_H}\ge \E[\|X-T_H(P)\|^2]\\-\E\left[\left(\|X-T_H(P)\|^2+\beta^2\right)^{q}\right]^{1/q}\P\left(\|X-T_H(P)\|>\beta\right)^{1-1/q}.
\end{multline}
\label{item:7:1}
\onepf{
We have
\begin{align*}
V_{\psi_H}&=\E\left[\psi_H\left(\|X-T_H(P)\|\right)^2\right]\\
&=\E[\beta^2 \wedge (\|X-T_H(P)\|)^2] \\
&=\E[\|X-T_H(P)\|^2]-\E[\left(\|X-T_H(P)\|^2+\beta^2\right)\1\{\|X-T_H(P)\|>\beta\}]
\end{align*}
Then, by Hölder inequality, we get the result announced.
}
\item We have $V_{\psi_H}\ge \E[\|X-\E[X]\|^2]-4^{q}\frac{\E\left[\|X-T_H(P)\|^{2q}\right]^{1-1/q}}{\left(\E\left[\|X-T_H(P)\|^{2q}\right]+\beta^{2q}\right)^{1-2/q}}.
$.\\
\onepf{
We use the following lemma
\begin{Lemma}~\label{lem:cantelli}
Let $Y$ be a positive real random variable, $\E[Y^q]<\infty$. We have for all $\lambda>0$,
$$
\P\left(Y\ge \lambda \right)\le 2^{q-1}\frac{\E\left[Y^q\right]}{\lambda^q+\E[Y^q]}
$$
\end{Lemma}
See Section~\ref{sec:proof_cantelli} for the proof. Then, for $Y=\|X-T_H(P)\|$ and from \ref{item:7:1}, we get,
\begin{multline*}
V_{\psi_H}\ge\E[\|X-T_H(P)\|^2]\\
-\E\left[\left(\|X-T_H(P)\|^2+\beta^2\right)^{q}\right]^{1/q}\left( 2^{2q-1}\frac{\E\left[\|X-T_H(P)\|^{2q}\right]}{\E\left[\|X-T_H(P)\|^{2q}\right]+\beta^{2q}} \right)^{1-1/q}.
\end{multline*}
Use the fact that $(a+b)^q\le 2^{q-1}(a^q+b^q)$,
\begin{align*}
V_{\psi_H}\ge&\E[\|X-T_H(P)\|^2]\\
&-2^{(q-1)/q+(2q-1)(1-1/q)}\frac{\E\left[\|X-T_H(P)\|^{2q}\right]^{1-1/q}}{\left(\E\left[\|X-T_H(P)\|^{2q}\right]+\beta^{2q}\right)^{1-2/q}}\\
\ge& \E[\|X-T_H(P)\|^2]-2^{2q}\frac{\E\left[\|X-T_H(P)\|^{2q}\right]^{1-1/q}}{\left(\E\left[\|X-T_H(P)\|^{2q}\right]+\beta^{2q}\right)^{1-2/q}}
\end{align*}
And finally, because $\E[X]$ is the minimizer of the quadratic loss,
\begin{align*}
V_{\psi_H}\ge \E[\|X-\E[X]\|^2]-4^{q}\frac{\E\left[\|X-T_H(P)\|^{2q}\right]^{1-1/q}}{\left(\E\left[\|X-T_H(P)\|^{2q}\right]+\beta^{2q}\right)^{1-2/q}}.
\end{align*}}
\item Similarly, if $X$ has $q>2$ finite moments, we have\\ $ v_{\psi_H}\ge \|\Sigma\|_{op} -4^{q}\frac{\E\left[\|X-T_H(P)\|^{2q}\right]^{1-1/q}}{\left(\E\left[\|X-T_H(P)\|^{2q}\right]+\beta^{2q}\right)^{1-2/q}}.$\\
\onepf{
Then, we operate the same manner for the bound on $v_{\psi_H}$. We have,
\begin{align*}
v_{\psi_H}=&\sup_{u \in S}\E\left[\frac{\langle u,X-T_H(P)\rangle ^2}{\|X-T_H(P)\|^2}\left(\beta^2 \wedge (\|X-T_H(P)\|)^2 \right)\right]\\
=&\sup_{u \in S}\left(\E\left[\langle u,X-T_H(P)\rangle ^2\right]\right.\\
&\left.-\E\left[\frac{\langle u,X-T_H(P)\rangle ^2}{\|X-T_H(P)\|^2}\left(\beta^2-(\|X-T_H(P)\|)^2 \right)\1\{\|X-T_H(P)\|\ge \beta \}\right]\right)
\end{align*}
Then, use Cauchy-Schwarz inequality,
\begin{align*}
v_{\psi_H}\ge& \sup_{u \in S}\E\left[\langle u,X-T_H(P)\rangle ^2\right]\\
&-\E\left[\left(\beta^2-(\|X-T_H(P)\|)^2 \right)\1\{\|X-T_H(P)\|\ge \beta \}\right]\\
\ge& \sup_{u \in S}\left(\E\left[\langle u,X-\E[X]\rangle ^2\right]+\langle u,\E[X]-T_H(P)\rangle^2\right.\\
&-\left. -\E\left[\left(\beta^2-(\|X-T_H(P)\|)^2 \right)\1\{\|X-T_H(P)\|\ge \beta \}\right]\right)\\
\ge& \sup_{u \in S}\E\left[\langle u,X-\E[X]\rangle ^2\right]\\
&-\E\left[\left(\beta^2-(\|X-T_H(P)\|)^2 \right)\1\{\|X-T_H(P)\|\ge \beta \}\right]
\end{align*}
Then, use the same reasoning as for the bound on $V_{\psi_H}$ to conclude that
$$ v_{\psi_H}\ge \|\Sigma\|_{op} -4^{q}\frac{\E\left[\|X-T_H(P)\|^{2q}\right]^{1-1/q}}{\left(\E\left[\|X-T_H(P)\|^{2q}\right]+\beta^{2q}\right)^{1-2/q}}$$
}
\end{npf}
\subsection{Proof of Lemma~\ref{lem:subgaussian}}\label{sec:proof_subgaussian}
From \cite[Theorem 4]{adamczak2008} and because $\|Y\|=\sup_{\|u\|=1}\langle Y,u\rangle $, there exists an absolute constant $C_1$ such that, for all $t\ge 0$,
\begin{multline*}
\P\left(\left\|\sum_{i=1}^n Y_i\right\|\ge \frac{3}{2}\E\left[\left\|\sum_{i=1}^n Y_i\right\|\right]+t \right)\\\le \exp\left( -\frac{t^2}{4n\sigma^2}\right)+3\exp\left(-\frac{t}{C_1\| \max_{1\le i\le n} \|Y_i\| \|_{\psi_1}} \right). 
\end{multline*}
where $\sigma^2=n\sup_{u \in S} \E[\langle X,u\rangle ^2]$.
Remark that $\sigma^2$ can be rewritten
$$
\sigma^2=n\sup_{u \in S} \langle u\E[X \otimes X], u\rangle =n\|\Sigma\|_{op}.
$$
Then, by Cauchy-Schwarz inequality,
$$\E\left[\left\|\sum_{i=1}^n Y_i\right\|\right]\le\E\left[\left\|\sum_{i=1}^n Y_i\right\|^2\right]^{1/2}=  \sqrt{n}\E\left[\left\|Y\right\|^2\right]^{1/2}.$$

\subsection{Proof of Lemma~\ref{lem:rough}}\label{sec:proof_rough}

$\rho$ is convex because $\rho''=\psi'\ge 0$ and it is increasing because $\psi=\rho'\ge 0$ ($\psi(0)=0$ and $\psi$ increasing). Then, from triangular inequality and Jensen's inequality, we have
$$ \rho\left(\frac{\|\E[X]-T(P)\|}{\beta}\right)\le \rho\left(\frac{\E[\|X-T(P)\|]}{\beta}\right)\le \E\left[ \rho\left(\frac{\|X-T(P)\|}{\beta}\right)\right].$$
By definition of $T(P)$, it is a minimizer of $\theta \mapsto \E\left[ \rho\left(\frac{\|X-\theta\|}{\beta}\right)\right]$, hence,
$$ \rho\left(\frac{\|\E[X]-T(P)\|}{\beta}\right)\le \E\left[ \rho\left(\frac{\|X-\E[X]\|}{\beta}\right)\right]$$
then, use the hypothesis to upper bound the right-hand side by $\rho(1/3)$, we get
$$ \rho\left(\frac{\|\E[X]-T(P)\|}{\beta}\right)\le \rho(1/3).$$
Finally, because $\rho$ is non-decreasing on $\R_+$ (its derivative is non-negative), we get the result.
%

\subsection{Proof of Lemma~\ref{lem:initialization}}\label{sec:proof_lem_initialization}
First, let us begin with $d=1$. We have that for all $t>0$,
\begin{align}\label{eq:med_init}
\P\left(Med(X_1^n) - \E_P[X] > t \right)&\le \P\left(\sum_{i=1}^n \1\{X_i - \E_P[X]> t\}\ge \frac{n}{2} \right)\nonumber\\
&\le  \P\left(\sum_{i=1}^n \1\{X'_i - \E_P[X]> t\}\ge \frac{n}{2} - |\mathcal{O}| \right)
\end{align}
By Hoeffding's inequality, we have
\begin{multline*}
\P\left(\sum_{i=1}^n  \1\{X_i - \E_P[X]> t\}\ge \frac{n}{2}-|\mathcal{O}|\right) \\
\le \exp\left(-2n\left(\frac{1}{2}-\frac{|\mathcal{O}|}{n}-\P\left(X - \E_P[X] > t \right) \right)^2  \right)
\end{multline*}
and by Chebychev inequality, for the choice $t=2\sqrt{2}\sigma$ we have \\$\P\left(X - \E_P[X] > t \right)\le 1/8$. Then, from Equation~\eqref{eq:med_init},
\begin{align*}
\P\left(Med(X_1^n) - \E_P[X] > 2\sqrt{2}\sigma \right) \le \exp\left(-2n\left(\frac{1}{2}-\frac{|\mathcal{O}|}{n}-\frac{1}{8} \right)^2\right) .
\end{align*}
and then, because $|\mathcal{O}|\le n/8$,
\begin{align*}
\P\left(Med(X_1^n) - \E_P[X] > 2\sqrt{2}\sigma \right) \le \exp\left(-n/8\right) .
\end{align*}
Now for dimension $d$, we use the dimension $1$ result on each coordinate, and by union bound we have that with probability larger than $1-de^{-n/8}$, for all $1\le j\le d$,
$$Med(\langle X_1^n, e_j\rangle)-\E_P[\langle X, e_j\rangle ]\le 2\sqrt{2}\sigma_j $$
and then, by taking the sum of the squares, we get $\|\theta_0-\E_P[X]\|^2\le 8\sum_{j=1}^d\sigma_j^2=8Tr(\Sigma).$
We conclude using that  $\|T(X_1^n)-\E_P[X]\|\le r_n$ with probability larger than $1-\delta$.

\subsection{Proof of Lemma~\ref{lem:surrogate}}\label{sec:proof_lem_surrogate}
The proof is derived from the proof of iterative reweighting algorithm for regression found in~\cite[Section 7.8]{robuststat}.

\textbf{First point}.\\
We have
\begin{multline*}
\textbf{U}_{\theta^{(m)}}(\theta)\\=\frac{1}{n}\sum_{i=1}^n\left( \frac{w_i(\theta^{(m)})}{2}\left(\frac{\|X_i-\theta\|}{\beta} \right)^2+ \rho_1(r_i(\theta^{(m)}))-\frac{1}{2}r_i(\theta^{(m)})\psi_1(r_i(\theta^{(m)}))\right). 
\end{multline*}
$\textbf{U}_{\theta^{(m)}}$ is a convex function, let us take its gradient to find its minimum,
$$\nabla \textbf{U}_{\theta^{(m)}}(\theta)=\frac{1}{n}\sum_{i=1}^n w_i(\theta^{(m)})\frac{\theta-X_i}{\beta^2} $$
Hence, the minimum is found for $\theta = \sum_{i=1}^n \frac{w_i(\theta^{(m)})}{\sum_{j=1}^n w_j(\theta^{(m)})} X_i=\theta^{(m+1)}$.

\textbf{Second point}.\\
For all $i\in \{1,\dots,n\}$,
$$g_i(x)=\frac{w_i(\kappa)}{2}x^2+ \rho_1(r_i(\kappa))-\frac{1}{2}r_i(\kappa)\psi_1(\kappa) - \rho_1(x)$$
We have that $g_i$ is differentiable and
$g_i'(x)=w_i(\kappa)x- \psi_1(x)=\frac{\psi_1(r_i(\kappa))}{r_i(\kappa)}x - \psi_1(x)$.
Then, having that $x \mapsto \psi_1(x)/x$ is non-increasing on $\R_+$, we have that for $x \in [0, r_i(\kappa)]$, $g_i'(x)\ge 0$ and for $x \ge r_i(\kappa)$, $g_i'(x)\le 0$. Hence, $g_i$ is minimal in $x=r_i(\kappa)$ and $g_i(x)\ge g_i(r_i(\kappa))=0$.

This prove that $g_i$ is a majorant of $\rho_1$ and by taking the sum, this implies that $\textbf{U}_\kappa$ is a majorant of $J_n$.

\textbf{Third point}.\\
Define
$$h_{i,\kappa}(\theta)= \frac{w_i(\kappa)}{2}\left(\frac{\|X_i-\theta\|}{\beta} \right)^2+ \rho_1(r_i(\kappa))-\frac{1}{2}r_i(\kappa)\psi_1(r_i(\kappa)) - \rho_1\left(\frac{\|X_i-\theta\|}{\beta}\right),$$
we have by definition of $w_i(\kappa)$,
\begin{equation}\label{eq:h0}
h_{i,\kappa}(\kappa)=\frac{w_i(\kappa)}{2}r_i(\kappa)^2+ \rho_1(r_i(\kappa))-\frac{1}{2}r_i(\kappa)\psi_1(r_i(\kappa)) - \rho_1(r_i(\kappa))=0
\end{equation}
and moreover, $h_i$ is a differentiable function whose gradient is
$$\nabla h_{i,\kappa}(\theta)=w_i(\kappa)\frac{(\theta-X_i)}{\beta^2} - \frac{\theta-X_i}{\beta\|\theta-X_i\|}\psi_1\left(\frac{\|X_i-\theta\|}{\beta}\right),$$
and we can verify that
\begin{equation}\label{eq:Dh0}
\nabla h_{i,\kappa}(\kappa)=w_i(\kappa)\frac{(\kappa-X_i)}{\beta^2} - \frac{\kappa-X_i}{\beta^2r_i}\psi\left(r_i\right)=0.
\end{equation}
Let us show that $\nabla h_\kappa$ is Lipshitz. We have

$$Hess(h_{i,\kappa})=w_i(\kappa)\frac{I_d}{\beta^2}-Hess(\theta \mapsto \rho_1(\|X_i-\theta\|))$$
Then, use that $\theta \mapsto \rho_1(\|X_i-\theta\|)$ is convex, hence its Hessian is positive and we have for all $u \in S$,
$$u^THess(h_{i,\kappa})u\le \frac{w_i}{\beta^2} $$
This conclude that $\nabla h_{i,\kappa}$ is $w_i/\beta^2$-Lipshitz continuous and hence by summing over $i$, $\nabla h_\kappa$ is Lipshitz continuous with Lipshitz constant $L=\frac{1}{n\beta^2}\sum_{i=1}^n \frac{\psi_1(r_i)}{r_i}\le 1/\beta^2$.

\textbf{Fourth point} already verified using Equations~\eqref{eq:h0} and \eqref{eq:Dh0}.

\subsection{Proof of Lemma~\ref{lem:descent}}\label{sec:proof_lem_descent}
The function $\textbf{U}_\kappa$ can be rewritten as follows:
\begin{align*}\textbf{U}_\kappa(\theta)=&\frac{1}{n}\sum_{i=1}^n\left( \frac{w_i(\kappa)}{2}\left(\frac{\|X_i-\theta\|}{\beta} \right)^2+ \rho_1(r_i(\kappa))-\frac{1}{2}r_i(\kappa)\psi_1(r_i(\kappa))\right)\\
=& \frac{1}{n}\sum_{i=1}^n \frac{w_i(\kappa)}{2}\frac{\|X_i-\kappa\|^2+\|\kappa-\theta\|+2\langle X_i-\kappa, \kappa-\theta\rangle }{\beta^2} \\
&+ \frac{1}{n}\sum_{i=1}^n\left( \rho_1(r_i(\kappa))-\frac{1}{2}r_i(\kappa)\psi_1(r_i(\kappa))\right)\\
=& \frac{1}{n}\sum_{i=1}^n \frac{r_i(\kappa)\psi_1(r_i(\kappa))}{2}+\frac{w_i(\kappa)}{2}\frac{\|\kappa-\theta\|+2\langle X_i-\kappa, \kappa-\theta\rangle }{\beta^2} \\
&+\frac{1}{n}\sum_{i=1}^n \left( \rho_1(r_i(\kappa))-\frac{1}{2}r_i(\kappa)\psi_1(r_i(\kappa))\right)\\
=& \frac{1}{n}\sum_{i=1}^n\left(\rho_1(r_i(\kappa)) +\frac{\langle X_i-\kappa, \kappa-\theta\rangle }{\beta\|X_i-\kappa\|}\psi_1(r_i(\kappa)) + \frac{w_i(\kappa)\|\kappa-\theta\|}{2\beta^2}\right)\\
=& J(\kappa)+2\langle \nabla J(\kappa), \theta - \kappa\rangle + \frac{1}{2n\beta^2}\sum_{i=1}^n w_i(\kappa)\|\theta -\kappa\|^2
\end{align*}
Then, let
$f(\theta)=\sum_{i=1}^n\frac{ w_i(\theta)}{\sum_{j=1}^n w_j(\theta)} X_i$ where $w_i(\theta)=\psi_1(r_i(\theta))$ and $r_i(\theta)=\|X_i-\theta\|/\beta$. Inject $\theta= f(\kappa)$ to get,
$$\textbf{U}_\kappa(f(\kappa)) =J(\kappa)+2\langle \nabla J(\kappa), f(\kappa) - \kappa\rangle + \frac{1}{2n\beta^2}\sum_{i=1}^n w_i(\kappa)\|f(\kappa) -\kappa\|^2  $$
Then, use that $ \textbf{U}_\kappa$ is a majorant of $J$ (Lemma~\ref{lem:surrogate}) and apply this to $\kappa = \theta^{(m)}$ to obtain
\begin{multline*}J(\theta^{(m+1)})\le J(\theta^{(m)})\\+\langle \nabla J(\theta^{(m)}), \theta^{(m+1)}-\theta^{(m)}\rangle + \frac{\sum_{j=1}^n w_j(\theta^{(m)})}{2n\beta^2}\|\theta^{(m+1)}-\theta^{(m)}\|^2 \end{multline*}
Then, because $J_n$ is convex, we have that for all $\theta \in \R^d$, $J(\theta^{(m)}\le J(\theta)+ \langle \nabla J(\theta^{(m)}), \theta^{(m)}-\theta \rangle $, hence
\begin{align*}
J(&\theta^{(m+1)})\\
\le& J(\theta)+\langle \nabla J(\theta^{(m)}), \theta^{(m+1)}-\theta\rangle + \frac{\sum_{j=1}^n w_j(\theta^{(m)})}{2n\beta^2}\|\theta^{(m+1)}-\theta^{(m)}\|^2 \\
=& J(\theta)+\langle \nabla U_{\theta^{(m)}}(\theta^{(m)}), \theta^{(m+1)}-\theta\rangle + \frac{\sum_{j=1}^n w_j(\theta^{(m)})}{2n\beta^2}\|\theta^{(m+1)}-\theta^{(m)}\|^2 \\
=& J(\theta)+\left\langle \frac{1}{n}\sum_{i=1}^n w_i(\theta^{(m)})(X_i - \theta^{(m)}) , \theta^{(m+1)}-\theta\right\rangle \\
&+ \frac{\sum_{j=1}^n w_j(\theta^{(m)})}{2n\beta^2}\|\theta^{(m+1)}-\theta^{(m)}\|^2 \\
=& J(\theta)+\sum_{j=1}^n w_j(\theta^{(m)}) \left\langle \theta^{(m+1)}-\theta^{(m)} , \theta^{(m+1)}-\theta\right\rangle\\
& + \frac{\sum_{j=1}^n w_j(\theta^{(m)})}{2n\beta^2}\|\theta^{(m+1)}-\theta^{(m)}\|^2
\end{align*}
Then, use the following identity: for all $u, v, w \in \R^d$,
$$2 \langle w - v , v - u\rangle = \|w-v\|^2  - \|w-u\|^2 + \|u-v\|^2 $$
to conclude that
$$J(\theta^{(m+1)}) \le  J(\theta)+\frac{1}{2n\beta^2}\sum_{j=1}^n w_j(\theta^{(m)}) \left(\|\theta^{(m)} - \theta\| - \|\theta^{(m+1)}-\theta\|\right) $$
\subsection{Proof of Lemma~\ref{lem:convex2}}\label{sec:proof_convex2}

By Assumption~\ref{ass:3}, we have
$$\frac{1}{n}\sum_{i=1}^n \psi_1'\left(\frac{\|X_i-\theta\|}{\beta}\right) \ge \frac{\gamma}{n}\sum_{i=1}^n \1\{\|X_i-\theta\|\le \beta\} $$
Then, because $\theta \in \Theta$ and $\beta \ge 2\sqrt{2 Tr(\Sigma)}+2r_n + \psi_1^{-1}\left(\sqrt{2V_{\psi_1}}\right)$,
\begin{align*}
\frac{1}{n}\sum_{i=1}^n \psi_1'\left(\frac{\|X_i-\theta\|}{\beta}\right) &\ge \frac{\gamma}{n}\sum_{i=1}^n \1\{\|X_i-T(X_1^n)\|\le \beta - 2\sqrt{2 Tr(\Sigma)} - r_n\} \\
&\ge \frac{\gamma}{n}\sum_{i=1}^n \1\{\|X_i-T(X_1^n)\|\le \psi_1^{-1}\left(\sqrt{2V_{\psi_1}}\right)+r_n\}
\end{align*}
Now, having that with probability larger than $1-\delta$, $\|T(X_1^n)-T(P)\| \le r_n $, we have with probability larger than $1-\delta$,
\begin{align*}
\frac{1}{n}\sum_{i=1}^n \psi_1'\left(\frac{\|X_i-\theta\|}{\beta}\right) &\ge \frac{\gamma}{n}\sum_{i=1}^n \1\{\|X_i-T(P)\|\le \psi_1^{-1}\left(\sqrt{2V_{\psi_1}}\right)\}.
\end{align*}
Then, as we don't have informations on the outliers, we take them out:
\begin{align*}
\frac{1}{n}\sum_{i=1}^n \psi_1'\left(\frac{\|X_i-\theta\|}{\beta}\right) &\ge \frac{\gamma}{n} \sum_{i =1}^n \1\left\{\|X_i'-T(P)\|\le \psi_1^{-1}\left(\sqrt{2V_{\psi_1}}\right)\right\}  - \gamma\frac{|\mathcal{O}|}{n}
\end{align*}

Then, by Hoeffding's inequality, we have that for $t>0$, with probability larger than $1-\delta-\exp(-2nt^2)$,
\begin{align*}
\frac{1}{n}\sum_{i=1}^n \psi_1'\left(\frac{\|X_i-\theta\|}{\beta}\right) &\ge \gamma  \P\left(\|X - T(P)\|\le \psi_1^{-1}\left(\sqrt{2V_{\psi_1}}\right) \right) -  \gamma t- \gamma\frac{|\mathcal{O}|}{n} \\
&\ge \gamma \P\left(\psi_1(\|X - T(P)\|)^2\le 2V_{\psi_1} \right) - \gamma t -\gamma \frac{|\mathcal{O}|}{n}
\end{align*}
By Markov inequality and having $|\mathcal{O}|\le n/8$, we have that
\begin{align*}
\frac{1}{n}\sum_{i=1}^n \psi_1'\left(\frac{\|X_i-\theta\|}{\beta}\right) &\ge \gamma \left(\frac{1}{2}-t-\frac{1}{8} \right).
\end{align*}
Hence, with probability larger than $1-\delta-e^{-n/32}$, we have
$$\frac{1}{n}\sum_{i=1}^n \psi_1'\left(\frac{\|X_i-\theta\|}{\beta}\right) \ge\frac{\gamma}{4} .
$$
\subsection{Proof of Lemma~\ref{lem:catoni_moments}}\label{sec:proof_catoni_moments}
For $q\in \N^*$, let
$$
g_q:x\mapsto \begin{cases}
q^qx/(e^{q}-1) & if\, x\in [0,e^{q}-1]\\
\log(1+x)^{q}& if\, x>e^{q}-1
\end{cases}$$

\begin{npf}
\item $g_q$ is a concave function over $\R_+$.\label{item:cat:1}\\
\onepf{
$g_q$ is continuous at $e^{q}-1$, the left and right limits are equal to $q^q$. $g_q$ is derivable on $[0,e^{q}-1)$ and $(e^{q}-1,\infty)$.
This derivative is non-increasing on both intervals.
At $e^{q}-1$, the left derivative is $q^q(e^{q}-1)^{-1}$ while the derivative on the right is $q^qe^{-q}$.
Thus, the left derivative at $e^{q}-1$ is larger than the right derivative.
Hence the derivative is non-increasing on $\R_+$, $g_q$ is concave on $\R_+$.}
\item $\E[\log(1+Z)^q]\le \E[g_q(Z)]\le g_q(\E[Z]) .$\label{item:cat:2}\\
\onepf{
By concavity, $\log(1+x)^q$ is smaller than its tangent in $e^{q-1}-1$.
This tangent is given by the function
$$x\mapsto (q-1)^q+\frac{q^{q}}{e^{q-1}-1}\left(x-(e^{q-1}-1) \right).$$
This last function is clearly smaller than the function $x\mapsto q^qx/(e^{q}-1)$.
Hence, $x\mapsto \log(1+x)^q$ is smaller than $g_q$, we found a concave upper bound of $x\mapsto \log(1+x)^q$.

Since $g_q$ is concave (\ref{item:cat:1}), by Jensen's inequality, for any positive random variable $Z$ such that $\E[Z]<\infty$, we have
$\E[\log(1+Z)^q]\le \E[g_q(Z)]\le g_q(\E[Z]) .$}
\item For any $q \in \N^*$, $\E[\psi_C(\|X-T_C(P)\|)^q]\le q!(\beta s)^q.
$\\
\onepf{Then, for all $x$, we have $g_q(x)\le \max(q^q,\log(1+x)^q)$, hence, by \ref{item:cat:2}
$$ \E[\log(1+Z)^q]\le \max(q^q,\log(1+\E[Z])^q).$$
Finally, use that $q^q\le q!e^q$ to get
\begin{equation}\label{eq:jensen_Z}
\E[\log(1+Z)^q]\le q!\max(e,\log(1+\E[Z]))^q.
\end{equation}
Denote $s=\max\left(e,\log\left(1+\frac{\E[\|X-T_C(P)\|])}{\beta}+\frac{\E[\|X-T_C(P)\|^2]}{2\beta^2} \right)\right),$
and apply equation~\eqref{eq:jensen_Z} to $Z=X/\beta+X^2/(2\beta^2)$ to get
$\E[\psi_C(\|X-T_C(P)\|)^q]\le q!(\beta s)^q.
$}
\end{npf}

\subsection{Proof of Lemma~\ref{lem:concentration_poly}}\label{sec:proof_concentration_poly}
By Markov's inequality, we have for any $\lambda>0$,
\begin{equation}\label{eq:poly_markov}
t_{\IF}(\lambda; T_P(P), X_1^n)\le \frac{\E\left[\left\|\frac{1}{n}\sum_{i=1}^n \frac{X_i-T_P(P)}{\|X_i-T_P(P)\|}\psi_P(\|X_i-T_P(P)\|)\right\|^{qp}\right]}{\lambda ^{qp}}.
\end{equation}
Let $Y_i=\frac{1}{n}\sum_{i=1}^n \frac{X_i-T_P(P)}{\|X_i-T_P(P)\|}\psi_P(\|X_i-T_P(P)\|)$ for $1\le i\le d$.
from \cite[Theorem 1.2.5]{Gine}, there exists an absolute constant $K>0$ such that
\begin{equation}\label{eq:gine_moments}
\E\left[\left\|\sum_{i=1} ^n Y_i \right\|^{pq}\right]^{1/(pq)}\le Kpq\left(\E\left[\left\|\sum_{i=1} ^n Y_i \right\|^{2}\right]^{1/2}+\E\left[\max_{1\le i\le n}\left\| Y_i \right\|^{pq}\right]^{1/(pq)} \right).
\end{equation}
 \begin{npf}
 \item We have $\E\left[\left\|\sum_{i=1} ^n Y_i \right\|^{2}\right]\le 4n\E[\|Y\|^{pq}]^{2/(pq)}.
$.\label{item:poly:1}\\
 \onepf{
Let $\varepsilon_1,\dots,\varepsilon_n$ denote i.i.d Rademacher random variable independents from $Y_1,\dots,Y_n$.
By the symmetrization lemma (see \cite[Lemma 1.2.6]{Gine}),
\begin{align*}
\E\left[\left\|\sum_{i=1} ^n Y_i \right\|^{2}\right]&\le 4\E\left[\left\|\sum_{i=1} ^n \varepsilon_i Y_i \right\|^{2}\right]= 4\E\left[\sum_{i=1} ^n \left\| Y_i \right\|^{2}\right]=4n\E[\|Y\|^2].
\end{align*}
Thus, by Jensen's inequality,
$
\E\left[\left\|\sum_{i=1} ^n Y_i \right\|^{2}\right]\le 4n\E[\|Y\|^{pq}]^{2/(pq)}.$}
\item $\E\left[\max_{1\le i\le n}\left\| Y_i \right\|^{pq}\right] \le n^{pq/2}\E\left[\left\| Y\right\|^{pq}\right]$. \label{item:poly:2}
\onepf{As the max of $n$ non-negative real numbers is smaller than their sum, we have
$$\E\left[\max_{1\le i\le n}\left\| Y_i \right\|^{pq}\right]\le  \E\left[\sum_{1\le i\le n}\left\| Y_i \right\|^{pq}\right]\le n\E\left[\left\| Y\right\|^{pq}\right]\le n^{pq/2}\E\left[\left\| Y\right\|^{pq}\right].$$}
\item For any $\lambda >0$, $t_{\IF}(\lambda; T_P(P),X_1^n)\le \frac{\E[\|X-T_P(P)\|^q]}{\beta^q}\left(\frac{Kpq\beta}{\sqrt{n}\lambda}\right)^{qp}$.\\
\onepf{From \ref{item:poly:1} and \ref{item:poly:2} and equation~\eqref{eq:gine_moments}, we get
\begin{equation}\label{eq:moments_poly}
\E\left[\left\|\sum_{i=1} ^n Y_i \right\|^{pq}\right]^{1/(pq)}\le 3Kpq\sqrt{n}\E\left[\left\| Y\right\|^{pq}\right]^{1/(pq)}.
\end{equation}
From equations~\eqref{eq:poly_markov} and~\eqref{eq:moments_poly} and if we re-inject the definition of $Y_i$'s, we get
$$t_{\IF}(\lambda; T_P(P),X_1^n)\le \E[\psi_P(\|X-T_P(P)\|)^{pq}]\left(\frac{Kpq}{\sqrt{n}\lambda}\right)^{qp}.$$
Then, use that $\psi_p(\|x\|)\le \|x\|^{1/p}\beta^{1-1/p}$ to get
$$t_{\IF}(\lambda; T_P(P),X_1^n)\le \frac{\E[\|X-T_P(P)\|^q]}{\beta^q}\left(\frac{Kpq\beta}{\sqrt{n}\lambda}\right)^{qp}.$$}
 \end{npf}

\subsection{Proof of Lemma~\ref{lem:cantelli}}\label{sec:proof_cantelli}
We have for all $u,\lambda>0$,
\begin{align*}
\P\left(Y\ge \lambda \right)&=\P\left((Y+u)^{q}\ge (\lambda+u)^{q} \right)\le \frac{\E\left[(Y+u)^{q}\right]}{(\lambda+u)^{q}}\le \frac{\E\left[(Y+u)^{q}\right]}{(\lambda^{q/2}+u^{q/2})^{2}}.
\end{align*}
Then, use that by convexity of the $q^{th}$-power function, $(a+b)^q \le 2^{q-1}(a^q+b^q)$ and also $(a+b)^q\ge a^q+b^q$,
\[
\P\left(Y\ge \lambda \right)\le 2^{q-1}\frac{\E\left[Y^q+u^q\right]}{(\lambda^{q/2}+u^{q/2})^{2}}.
\]
Take $u=\E[Y^{q}]^{2/q}/\lambda$ to get,
\begin{align*}
\P\left(Y\ge \lambda \right)&\le 2^{q-1}\frac{\E\left[Y^q\right]+\frac{\E[Y^q]^2}{\lambda^q}}{\lambda^q(1+\frac{\E[Y^q]}{\lambda^q})^{2}}=2^{q-1}\frac{\E\left[Y^q\right]}{\lambda^q(1+\frac{\E[Y^q]}{\lambda^q})}=2^{q-1}\frac{\E\left[Y^q\right]}{\lambda^q+\E[Y^q]}.
\end{align*}

\section{Technical tools}\label{sec:reminder}
We remind the reader of Bernstein inequality, a classical concentration inequality, this form of Bernstein inequality is borrowed from \cite[Theorem 2.10]{concentration}.
\begin{Theorem}
Let $X_1,\dots,X_n$ be independent real-valued random variables. Assume that there exist positive numbers $v$ and $c$ such that $\sum_{i=1}^n \E[X_i^2]\le v$ and
$$\sum_{i=1}^n \E[(X_i)_+^q]\le vc^{q-2} \quad \text{for all integers }q\ge 3,$$
where $x_+=\max(0,x)$.
Then for all $t>0$
$$\P\left(\sum_{i=1}^n (X_i-\E[X_i]) \ge \sqrt{2vt}+ct\right)\le e^{-t}.$$
\end{Theorem}

The following theorem is borrowed from \cite[ Theorem 4]{adamczak2008}, it is a concentration inequality for suprema of sums of independent random variables.
\begin{Theorem}
Let $X_1,\dots,X_n$ be independent random variables with values in a measurable space $(\mathcal{S},\mathcal{B})$ and let $\mathcal{F}$ be a countable class of measurable functions $f: \mathcal{S}\to \R$. Assume that for every $f\in \mathcal{F}$ and every $i$, $\E[f(X_i)]=0$ and for any $\alpha \in (0,1]$ and all $i$, $\|\sup_{f}|f(X_i)|\|_{\psi_\alpha}<\infty$. Let
$$Z=\sup_{f \in \mathcal{F}}\left|\sum_{i=1}^n f(X_i) \right|. $$
Define moreover
$$\sigma^2=\sup_{f \in \mathcal{F}}\sum_{i=1}^n \E[f(X_i)^2].$$
Then, for all $0<\eta<1$ and $\delta>0$, there exists a constant $C=C(\alpha,\eta,\delta)>0$ such that for all $t\ge 0$,
\begin{multline}
\P(Z\ge (1+\eta)\E[Z]+t)\le \\
\exp\left(-\frac{t^2}{2(1+\delta)\sigma^2} \right)+3\exp\left( -\left(\frac{t}{C\|\max_i \sup_{f\in \mathcal{F}}|f(X_i)|  \|_{\psi_\alpha}} \right)^{\alpha} \right),
\end{multline}
and
\begin{multline}
\P(Z\le (1-\eta)\E[Z]-t)\le \\
\exp\left(-\frac{t^2}{2(1+\delta)\sigma^2} \right)+3\exp\left( -\left(\frac{t}{C\|\max_i \sup_{f\in \mathcal{F}}|f(X_i)|  \|_{\psi_\alpha}} \right)^{\alpha} \right).
\end{multline}
\end{Theorem}

%
\end{appendix}

\section*{Acknowledgements}
The author would like to thank Elvezio Ronchetti, for the discussions that gave him the ideas behind this article, and his two phd-advisors Matthieu Lerasle and Guillaume Lecué for their advice and reviews.
%
%



\bibliographystyle{imsart-number} 
\bibliography{biblio_huber}       


\end{document}